\documentclass[12pt]{article}
\usepackage{amssymb, amsmath, amsthm, amscd}
\usepackage[utf8]{inputenc}
\usepackage[all,cmtip]{xy}
\usepackage{enumitem}
\usepackage{setspace}
\usepackage{mathdots}

\usepackage[mathscr]{euscript}
\usepackage[colorlinks=true]{hyperref}
\hypersetup{colorlinks   = true}
\hypersetup{linkcolor=blue}
\usepackage{tikz}

\begin{document}

\mathchardef\mhyphen="2D
\newtheorem{The}{Theorem}[section]
\newtheorem{Lem}[The]{Lemma}
\newtheorem{Prop}[The]{Proposition}
\newtheorem{Cor}[The]{Corollary}
\newtheorem{Rem}[The]{Remark}
\newtheorem{Obs}[The]{Observation}
\newtheorem{SConj}[The]{Standard Conjecture}
\newtheorem{Titre}[The]{\!\!\!\! }
\newtheorem{Conj}[The]{Conjecture}
\newtheorem{Question}[The]{Question}
\newtheorem{Prob}[The]{Problem}
\newtheorem{Def}[The]{Definition}
\newtheorem{Not}[The]{Notation}
\newtheorem{Claim}[The]{Claim}
\newtheorem{Conc}[The]{Conclusion}
\newtheorem{Ex}[The]{Example}
\newtheorem{Fact}[The]{Fact}
\newtheorem{Formula}[The]{Formula}
\newtheorem{Formulae}[The]{Formulae}
\newtheorem{The-Def}[The]{Theorem and Definition}
\newtheorem{Prop-Def}[The]{Proposition and Definition}
\newtheorem{Lem-Def}[The]{Lemma and Definition}
\newtheorem{Cor-Def}[The]{Corollary and Definition}
\newtheorem{Conc-Def}[The]{Conclusion and Definition}
\newtheorem{Terminology}[The]{Note on terminology}
\newcommand{\C}{\mathbb{C}}
\newcommand{\R}{\mathbb{R}}
\newcommand{\N}{\mathbb{N}}
\newcommand{\Z}{\mathbb{Z}}
\newcommand{\Q}{\mathbb{Q}}
\newcommand{\Proj}{\mathbb{P}}
\newcommand{\Rc}{\mathcal{R}}
\newcommand{\Oc}{\mathcal{O}}
\newcommand{\Vc}{\mathcal{V}}
\newcommand{\Id}{\operatorname{Id}}
\newcommand{\pr}{\operatorname{pr}}
\newcommand{\rk}{\operatorname{rk}}
\newcommand{\del}{\partial}
\newcommand{\delbar}{\bar{\partial}}
\newcommand{\Cdot}{{\raisebox{-0.7ex}[0pt][0pt]{\scalebox{2.0}{$\cdot$}}}}
\newcommand\nilm{\Gamma\backslash G}
\newcommand\frg{{\mathfrak g}}
\newcommand{\fg}{\mathfrak g}
\newcommand{\Oh}{\mathcal{O}}
\newcommand{\Kur}{\operatorname{Kur}}
\newcommand\gc{\frg_\mathbb{C}}
\newcommand\hisashi[1]{{\textcolor{red}{#1}}}
\newcommand\dan[1]{{\textcolor{blue}{#1}}}
\newcommand\luis[1]{{\textcolor{orange}{#1}}}

\begin{center}

{\Large\bf Generalised Hermite-Einstein Fibre Metrics and Slope Stability for Holomorphic Vector Bundles}

\end{center}

\begin{center}

{\large Dan Popovici}

\end{center}

\vspace{1ex}

\noindent{\small{\bf Abstract.} Let $X$ be a compact complex manifold of dimension $n$ and let $m$ be a positive integer with $m\leq n$. Assume that $X$ admits a K\"ahler metric $\omega$ and a weakly positive, $\partial\bar\partial$-closed, smooth $(n-m,\,n-m)$-form $\Omega$. We introduce the notions of $(\omega,\,\Omega)$-Hermite-Einstein holomorphic vector bundles and $(\omega,\,\Omega)$(-semi)-stable coherent sheaves on $X$ by generalising the classical definitions depending only on $\omega$. We then prove that the $(\omega,\,\Omega)$-Hermite-Einstein condition implies the $(\omega,\,\Omega)$-semi-stability of a holomorphic vector bundle and its splitting into $(\omega,\,\Omega)$-stable subbundles. This extends a classical result by Kobayashi and L\"ubke to our generalised setting. In the appendix, we propose notions of both strongly and weakly (strictly) positive forms and currents and discuss their various properties.}

\vspace{1ex}

\section{Introduction}\label{section:introduction} In our recent joint study with S. Dinew ([DP25a] and [DP25b]) of $m$-positivity, we gave a generalisation of Lamari's duality criterion to arbitrary bidegrees $(m,\,m)$. Specifically, let $X$ be a compact complex manifold with $\mbox{dim}_\C X = n$, let $m\in\{1,\dots , n\}$ and let $\theta\in C_{m,\,m}^\infty(X,\,\R)$. (By $C_{p,\,q}^\infty(X,\,\mathbb{K})$, we mean the space of $C^\infty$ differential forms of bidegree $(p,\,q)$ on $X$ with values in $\mathbb{K} = \R$ or $\C$.) Then, Lemma 2.1 of [DP25b] asserts that the following statements are equivalent:

\vspace{1ex}

(i)\, There exists a real current $S$ of bidegree $(m-1,\,m-1)$ on $X$ such that the $(m,\,m)$-current $\theta + i\partial\bar\partial S$ is strongly semi-positive (a fact denoted by $\theta + i\partial\bar\partial S\geq 0$ (strongly)) on $X$;

\vspace{1ex}

(ii)\, $\int\limits_X\theta\wedge\Omega \geq 0$ for every $\Omega\in C_{n-m,\,n-m}^\infty(X,\,\R)$ such that $\partial\bar\partial\Omega = 0$ and $\Omega>0$ (weakly) on $X$.

\vspace{1ex}

The weak (strict) positivity condition satisfied by $\Omega$ is to be construed in either of the equivalent senses introduced in the appendix ($\S$\ref{section:appendix}): {\bf metrically weak (strict) positivity} (cf. $(2)$ of Definition \ref{Def:weak-strict-positivity_pp-forms_v-space}) or {\bf decomposably weak (strict) positivity} (cf. Definition \ref{Def:decomposably-weak-strict-positivity_pp-forms_v-space}). The equivalence of these two properties is proved in Proposition \ref{Prop:weak-pos_nortions-hierarchy}, while the duality between strongly semi-positive currents of any bidegree $(q,\,q)$ and (metrically or decomposably) weakly positive smooth forms of the complementary bidegree $(n-q,\,n-q)$ is proved in Theorem \ref{The:duality_decomp-weakly_s-semi-pos-currents}.

\vspace{1ex}

In the present paper, we use the flexibility afforded by a test form $\Omega\in C_{n-m,\,n-m}^\infty(X,\,\R)$ with the above two properties to generalise two classical notions:

\vspace{1ex}

$(1)$\, Hermite-Einstein fibre metrics on a holomorphic vector bundle;

\vspace{1ex}

$(2)$\, slope stability of torsion-free coherent sheaves.

\vspace{1ex}

We will often assume the existence of a K\"ahler metric $\omega$ on $X$ so that certain integrals depend only on the cohomology classes (in the Bott-Chern or Aeppli senses) of the forms involved; this includes the {\bf $(\omega,\,\Omega)$-degree} of a holomorphic line bundle $L$ over $X$ introduced in Definition \ref{Def:stability-gen}.

Now, there is a huge literature on both the Hermite-Einstein theory and various notions of stability. In order to keep the exposition as focused as possible, we only cite the references that are used in our discussion and follow, by adapting it to our generalised setting, the presentation in the L\"ubke-Teleman book [LT95].

\vspace{1ex}

Having assumed the existence of a K\"ahler metric $\omega$ and of a form $\Omega\in C_{n-m,\,n-m}^\infty(X,\,\R)$, for some $m\in\{1,\dots , n\}$, such that $\partial\bar\partial\Omega = 0$ and $\Omega>0$ (metrically weakly) on $X$, we define:

\vspace{1ex}

$\bullet$ the notion of {\bf $(\omega,\,\Omega)$-Hermite-Einstein} $C^\infty$ fibre metric $h$ on a holomorphic vector bundle $E$ of rank $r\geq 1$ on $X$ (cf. Definition \ref{Def:H-E}); 

\vspace{1ex}

$\bullet$ the notion of {\bf $(\omega,\,\Omega)$(-semi)-stability} for a coherent sheaf ${\cal F}$ of ${\cal O}_X$-modules on $X$ (cf. Definition \ref{Def:stability-gen}).

\vspace{1ex}

Our main result (Theorem \ref{The:stability-H-E_implication_gen}) asserts that the classical Kobayashi-L\"ubke theorem of [Kob82] and [Lub83] continues to hold in our generalised context. More precisely, the existence of an $(\omega,\,\Omega)$-Hermite-Einstein $C^\infty$ fibre metric $h$ on a holomorphic vector bundle $E$ of rank $r\geq 1$ over $X$ implies the $(\omega,\,\Omega)$-semi-stability of $E$ and its splitting as a direct sum $E=\oplus_i E_i$ of $(\omega,\,\Omega)$-stable holomorphic subbundles $E_i\subset E$.

The Kobayashi-Hitchin correspondence and Hermitian Yang-Mills connections on possibly non-K\"ahler manifolds have been studied by several authors, including Li and Yau in [LY86] and Buchdahl in [Buc88]. The present paper is geared towards a future development of applications of $m$-positivity to stability issues.

In future work, we hope to study various examples of Hermite-Einstein and (semi-)stable vector bundles in the generalised sense introduced in this paper, as well as the associated moduli spaces and their dependence on the parameters $\Omega$ and $m$.

\vspace{1ex}

We included the appendix ($\S$\ref{section:appendix}) in order to fill an apparent gap in the literature concerning the (semi-)positivity terminology for both forms and curents. After briefly recalling Lelong's definitions of (strongly and weakly) semi-positive forms and currents as presented in [Dem97, chapter III], we introduce the notions of (strict) positivity in both the strong and weak senses and study their mutual relations. The resulting classification turns out to be richer than expected. For example, we identify two non-equivalent notions of weakly (strictly) positive forms at a given point (see Definition \ref{Def:weak-strict-positivity_pp-forms_v-space} and Proposition \ref{Prop:met-weakly-pos_implies_dual-weakly-pos}) and similarly for currents on a manifold (see Definition \ref{Def:strict-positivity_currents} and Observation \ref{Obs:hierarchy_strict-pos_currents}).

\vspace{2ex}

\noindent {\bf Acknowledgments.} The author is very grateful to the referees for their careful reading of the manuscript and their helpful suggestions. In particular, the appendix -- the most significant addition to the first version of this paper -- was written at the request of one of the referees. The author also thanks S. Dinew and L. Ugarte, collaborators on other recent papers, for helpful discussions on the content of the appendix. He is also grateful to B. Berndtsson for sharing his thoughts on the appendix and for pointing out the reference [Xia26].

\section{Definitions and results}\label{section:main-def}

Let $E\longrightarrow X$ be a holomorphic vector bundle of rank $r\geq 1$ on a compact complex $n$-dimensional manifold and let $m\in\{1,\dots , n\}$. Suppose there exist a K\"ahler metric $\omega$ on $X$ and a real form $\Omega\in C^\infty_{n-m,\,n-m}(X,\,\R)$ such that \begin{eqnarray}\label{eqn:hypotheses_Omega}\Omega>0 \hspace{2ex} \mbox{(metrically weakly) \hspace{3ex} and \hspace{3ex}} \partial\bar\partial\Omega = 0\end{eqnarray} at every point of $X$.

For any $C^\infty$ Hermitian fibre metric $h$ on $E$, $\Theta_h(E)\in C^\infty_{1,\,1}(X,\,\mbox{End}\,(E))$ will denote the curvature form of the Chern connection $D_h = D_h' + \bar\partial$ of $(E,\,h)$. Thus, $D_h^2 = \Theta_h(E)\wedge\cdot$.
Meanwhile, $dV_\omega=\omega_n:=\omega^n/n!\in C^\infty_{n,\,n}(X,\,\R)$ will denote the volume form induced by $\omega$ on $X$.


\vspace{2ex}

We will generalise to this setting the classical notions of Hermite-Einstein metrics and slope stability, as well as the relations between these two theories.

The following notation from [Dem97] will be used. If $(E,\,h)$ is a complex Hermitian vector bundle of rank $r$ on a complex manifold $X$, for any $E$-valued forms $s,\,t$ of respective degrees $k$ and $l$, $\{s,\,t\} = \{s,\,t\}_h$ stands for the scalar-valued form of degree $k+l$ obtained by combining the wedge product of scalar-valued forms with the pointwise inner product $\langle\,\cdot\,,\,\cdot\,\rangle = \langle\,\cdot\,,\,\cdot\,\rangle_h$ defined by $h$ on the fibres of $E$. Specifically, if $\{e_1,\dots , e_r\}$ is a local frame for $E$ and $s=\sum_\lambda s_\lambda\otimes e_\lambda$ and $t=\sum_\mu t_\mu\otimes e_\mu$ are the corresponding local expressions of $s,t$, with the $s_\lambda$ scalar-valued forms of degree $k$ and the $t_\mu$ scalar-valued forms of degree $l$, one defines \begin{eqnarray*}\{s,\,t\}_h = \sum\limits_{\lambda,\,\mu=1}^r s_\lambda\wedge\bar{t}_\mu\,\langle e_\lambda,\,e_\mu\rangle_h.\end{eqnarray*} In particular, when $k=l=0$ (i.e. $s$ and $t$ are sections of $E$), one has: $\{s,\,t\}_h = \langle\,s,\,t\rangle_h$.   

\subsection{Generalised Hermite-Einstein fibre metrics}\label{subsection:H-E_gen}

We propose the following generalisation of the classical notion of (weakly)
Hermite-Einstein metrics by letting $\Omega$ replace $\omega^{n-m}$. Its dependence on both $\omega$ and $\Omega$ gives extra flexibility to this definition.   

\begin{Def}\label{Def:H-E} A $C^\infty$ Hermitian fibre metric $h$ on $E$ is said to be:

\vspace{1ex}

(i)\, {\bf weakly $(\omega,\,\Omega)$-Hermite-Einstein} if there exists a $C^\infty$ function $\lambda_h:X\longrightarrow\R$, called the {\bf Einstein factor}, such that \begin{eqnarray}\label{eqn:H-E}i\Theta_h(E)\wedge\omega^{m-1}\wedge\Omega =  \lambda_h\,dV_\omega\otimes\mbox{Id}_E,\end{eqnarray} where $\mbox{Id}_E$ is the identity endomorphism of $E$.

\vspace{1ex}

(ii)\, {\bf $(\omega,\,\Omega)$-Hermite-Einstein} if there exists a constant $\lambda_h\in\R$, called the {\bf Einstein factor}, such that condition (\ref{eqn:H-E}) is satisfied.

\end{Def}

This definition prompts the introduction of the following Laplace-type differential operator depending on $\omega$ and $\Omega$ and acting on $\R$-valued or $\C$-valued $C^\infty$ functions on $X$. It first appeared in [DP25b].

\begin{Def}\label{Def:P_omega_Omega} (Definition 2.7. in [DP25b]) In the above context, we set: \begin{eqnarray*}   P = P_{\omega,\,\Omega}:C^\infty(X,\,\R)\longrightarrow C^\infty(X,\,\R), \hspace{5ex} P_{\omega,\,\Omega}\varphi:=-\frac{i\partial\bar\partial\varphi\wedge\omega^{m-1}\wedge\Omega}{dV_\omega}.\end{eqnarray*}

The same definition is made for $P = P_{\omega,\,\Omega}:C^\infty(X,\,\C)\longrightarrow C^\infty(X,\,\C)$.
  
\end{Def}

The adjoint operator is computed in

\begin{Lem}\label{Lem:P-operator_adjoint} (Lemma 2.8. in [DP25b]) Let $P^\star_{\omega,\,\Omega}:C^\infty(X,\,\C)\longrightarrow C^\infty(X,\,\C)$ be the $L^2_\omega$-adjoint of the operator introduced in Definition \ref{Def:P_omega_Omega}. For every $C^\infty$ function $\varphi:X\longrightarrow\C$, the following identity holds when $X$ is compact: \begin{eqnarray*}\bigg(P^\star_{\omega,\,\Omega} - P_{\omega,\,\Omega}\bigg)\,\varphi = \frac{i(\bar\partial\varphi\wedge\partial\Omega - \partial\varphi\wedge\bar\partial\Omega)\wedge\omega^{m-1}}{dV_\omega}.\end{eqnarray*}

In particular, $P^\star_{\omega,\,\Omega}$ and $P_{\omega,\,\Omega}$ differ by a first-order operator, so $P^\star_{\omega,\,\Omega}$ is elliptic of order two with no zero$^{th}$-order terms.

\end{Lem}

\noindent {\it Proof.} Let $\varphi,\psi:X\longrightarrow\C$ be arbitrary $C^\infty$ functions. We have: \begin{eqnarray*}\langle\langle P_{\omega,\,\Omega}\varphi,\,\psi\rangle\rangle = -\int\limits_X\overline\psi\, i\partial\bar\partial\varphi\wedge\omega^{m-1}\wedge\Omega  \hspace{3ex}\mbox{and}\hspace{3ex} \langle\langle \varphi,\,P_{\omega,\,\Omega}\psi\rangle\rangle = -\int\limits_X\varphi\, i\partial\bar\partial\,\overline\psi\wedge\omega^{m-1}\wedge\Omega.\end{eqnarray*}

On the other hand, we have: \begin{eqnarray*}i\partial\bar\partial(\varphi\,\overline\psi) = \varphi\,i\partial\bar\partial\,\overline\psi + i\partial\varphi\wedge\bar\partial\,\overline\psi + \overline\psi\,i\partial\bar\partial\varphi + i\partial\overline\psi\wedge\bar\partial\varphi.\end{eqnarray*}

Plugging the value for $\varphi\,i\partial\bar\partial\,\overline\psi$ given by this equality into one of the previous expressions, we get the second equality below: \begin{eqnarray}\label{eqn:P-operator-adjoint_proof_1}\nonumber\langle\langle P^\star_{\omega,\,\Omega}\varphi,\,\psi\rangle\rangle = \langle\langle \varphi,\,P_{\omega,\,\Omega}\psi\rangle\rangle & = & \int\limits_X \overline\psi\,i\partial\bar\partial\varphi\wedge\omega^{m-1}\wedge\Omega + \int\limits_X i\partial\varphi\wedge\bar\partial\,\overline\psi\wedge\omega^{m-1}\wedge\Omega \\
  & + & \int\limits_X i\partial\overline\psi\wedge\bar\partial\varphi\wedge\omega^{m-1}\wedge\Omega - \int\limits_X i\partial\bar\partial(\varphi\,\overline\psi)\wedge\omega^{m-1}\wedge\Omega.\end{eqnarray}

Now, thanks to $\omega^{m-1}$ being $\partial$-closed and $\bar\partial$-closed and to $\partial\bar\partial\Omega = 0$, two applications of the Stokes theorem yield the following result for the last integral: \begin{eqnarray*}\int\limits_X i\partial\bar\partial(\varphi\,\overline\psi)\wedge\omega^{m-1}\wedge\Omega & = & i\, \int\limits_X\partial\bigg(\bar\partial(\varphi\,\overline\psi)\wedge\omega^{m-1}\wedge\Omega\bigg) + i\, \int\limits_X\bar\partial(\varphi\,\overline\psi)\wedge\omega^{m-1}\wedge\partial\Omega \\
  & = & i\, \int\limits_X\bar\partial\bigg(\varphi\,\overline\psi\,\omega^{m-1}\wedge\partial\Omega\bigg) + i\, \int\limits_X\varphi\,\overline\psi\,\omega^{m-1}\wedge\partial\bar\partial\Omega = 0.\end{eqnarray*}

Similarly, for the last integral on the first line of (\ref{eqn:P-operator-adjoint_proof_1}), we get: \begin{eqnarray*}\int\limits_X i\partial\varphi\wedge\bar\partial\,\overline\psi\wedge\omega^{m-1}\wedge\Omega = -\int\limits_X\overline\psi\, i\partial\bar\partial\varphi\wedge\omega^{m-1}\wedge\Omega -i\,\int\limits_X\overline\psi\,\partial\varphi\wedge\omega^{m-1}\wedge\bar\partial\Omega,\end{eqnarray*} (so, this computes the sum of the last two integrals on the first line of (\ref{eqn:P-operator-adjoint_proof_1})), while for the first integral on the second line of (\ref{eqn:P-operator-adjoint_proof_1}), we get: \begin{eqnarray*}\int\limits_X i\partial\overline\psi\wedge\bar\partial\varphi\wedge\omega^{m-1}\wedge\Omega & = & -\int\limits_X\overline\psi\, i\partial\bar\partial\varphi\wedge\omega^{m-1}\wedge\Omega + i\,\int\limits_X\overline\psi\,\bar\partial\varphi\wedge\omega^{m-1}\wedge\partial\Omega \\
  & = & \langle\langle P_{\omega,\,\Omega}\varphi,\,\psi\rangle\rangle + i\,\int\limits_X\overline\psi\,\bar\partial\varphi\wedge\omega^{m-1}\wedge\partial\Omega.\end{eqnarray*}

All these computation results transform (\ref{eqn:P-operator-adjoint_proof_1}) into: \begin{eqnarray*}\langle\langle P^\star_{\omega,\,\Omega}\varphi,\,\psi\rangle\rangle = \langle\langle P_{\omega,\,\Omega}\varphi,\,\psi\rangle\rangle + \int\limits_X \overline\psi\,i(\bar\partial\varphi\wedge\partial\Omega - \partial\varphi\wedge\bar\partial\Omega)\wedge\omega^{m-1}.\end{eqnarray*} This equality holds for all $C^\infty$ functions $\varphi,\psi:X\longrightarrow\C$, proving the stated formula for $P^\star_{\omega,\,\Omega} - P_{\omega,\,\Omega}$.

Since $P^\star_{\omega,\,\Omega}$ and $P_{\omega,\,\Omega}$ are differential operators of order two and differ by a first-order operator, they have the same principal part and the same zero-th order terms. Since $P_{\omega,\,\Omega}$ is elliptic of order two with no zero-th order terms, the same holds for $P^\star_{\omega,\,\Omega}$.

This last property of $P_{\omega,\,\Omega}$ can be seen in the following way. By Proposition \ref{Prop:Omega-prod-omega_m-1}, the $(n-1,\,n-1)$-form $\omega^{m-1}\wedge\Omega$ is positive definite at every point of $X$. Now, a simple fact in linear algebra (often attributed to [Mic83]) ensures the existence of a unique $C^\infty$ $(1,\,1)$-form $\rho$ that is positive definite and satisfies $\omega^{m-1}\wedge\Omega = \rho_{n-1}: = \rho^{n-1}/(n-1)!$ at every point of $X$. Thus, we get: \begin{eqnarray*}P_{\omega,\,\Omega}\varphi:=-\frac{i\partial\bar\partial\varphi\wedge\rho_{n-1}}{dV_\rho}\,\frac{dV_\rho}{dV_\omega} = f_{\rho,\,\omega}\,\bigg(-\Lambda_\rho(i\partial\bar\partial\varphi)\bigg) = f_{\rho,\,\omega}\,\Delta_\rho\varphi,\end{eqnarray*} where $f_{\rho,\,\omega}:=dV_\rho/dV_\omega$ is a $C^\infty$ function on $X$ with positive values and $\Delta_\rho:=-\Lambda_\rho(i\partial\bar\partial)$ is the usual Laplacian on functions induced by the Hermitian metric $\rho$. As can be seen, $\Delta_\rho$ has no zero-th order terms and is known to be elliptic.    \hfill $\Box$

\begin{Cor}\label{Cor:2-space-decomposition} (Corollary 2.9. in [DP25b]) If the manifold $X$ is compact,

\vspace{1ex}  

\noindent $\ker\,\bigg(P_{\omega,\,\Omega}:C^\infty(X,\,\R)\longrightarrow C^\infty(X,\,\R\bigg) = \R$ \hspace{1ex} and \hspace{1ex} $\ker\,\bigg(P_{\omega,\,\Omega}:C^\infty(X,\,\C)\longrightarrow C^\infty(X,\,\C)\bigg) =\C$.

\vspace{1ex}

\noindent Moreover, the $L^2_\omega$-orthogonal decompositions hold: \begin{eqnarray}\label{eqn:2-space-decomposition}C^\infty(X,\,\R) = \R\oplus\mbox{Im}\,P_{\omega,\,\Omega} \hspace{3ex}\mbox{and}\hspace{3ex} C^\infty(X,\,\C) = \C\oplus\mbox{Im}\,P_{\omega,\,\Omega}.\end{eqnarray}

\end{Cor}  

\noindent {\it Proof.} The statement follows from Lemma \ref{Lem:P-operator_adjoint} and the maximum principle applied to $P^\star_{\omega,\,\Omega}$ on the compact $X$. \hfill $\Box$

\vspace{2ex}

As in the standard case (see e.g. Lemma 2.1.5 in [LT95]), we shall now see that the existence of a weakly $(\omega,\,\Omega)$-Hermite-Einstein metric $h$ on $E$ is actually equivalent to the existence of an $(\omega,\,\Omega)$-Hermite-Einstein metric $\tilde{h}$ on $E$.

\begin{Prop}\label{Prop:equiv_w-H-E_H-E} Consider the above setting of a holomorphic vector bundle $E\longrightarrow (X,\,\omega,\,\Omega)$ over a compact K\"ahler manifold equipped with an extra form $\Omega\in C^\infty_{n-m,\,n-m}(X,\,\R)$ satisfying (\ref{eqn:hypotheses_Omega}).

  If $h$ is a {\bf weakly $(\omega,\,\Omega)$-Hermite-Einstein} fibre metric on $E$, there exists a $C^\infty$ function $f:X\longrightarrow\R$, unique up to an additive real constant, such that the rescaled fibre metric $h\,e^{-f}$ is {\bf $(\omega,\,\Omega)$-Hermite-Einstein}.

\end{Prop}

\noindent {\it Proof.} By assumption, there exists a $C^\infty$ function $\lambda_h:X\longrightarrow\R$ such that \begin{eqnarray*}i\Theta_h(E)\wedge\omega^{m-1}\wedge\Omega =  \lambda_h\,dV_\omega\otimes\mbox{Id}_E.\end{eqnarray*}

On the one hand, the $L^2_\omega$-orthogonal splitting (\ref{eqn:2-space-decomposition}) yields a unique pair $(c,\,P_{\omega,\,\Omega}f)\in\R\oplus\mbox{Im}\,P_{\omega,\,\Omega}$ such that $\lambda_h = c + P_{\omega,\,\Omega}f$. Since the $2^{nd}$-order differential operator $P_{\omega,\,\Omega}$ is elliptic with no zero$^{th}$-order term and $X$ is compact, the $C^\infty$ function $f:X\longrightarrow\R$ whose image under $P_{\omega,\,\Omega}$ has been specified is unique up to an additive real constant. 

On the other hand, for every $C^\infty$ function $f:X\longrightarrow\R$ (in particular, for the one found above, unique up to an additive real constant), we have:  \begin{eqnarray*}i\Theta_{h\,e^{-f}}(E) = i\Theta_h(E) + (i\partial\bar\partial f)\otimes\mbox{Id}_E.\end{eqnarray*} Hence, multiplying by $\omega^{m-1}\wedge\Omega$, we get the first equality below: \begin{eqnarray*}i\Theta_{h\,e^{-f}}(E)\wedge\omega^{m-1}\wedge\Omega =
  (\lambda_h - P_{\omega,\,\Omega}f)\,dV_\omega\otimes\mbox{Id}_E = c\,dV_\omega\otimes\mbox{Id}_E .\end{eqnarray*}

This proves that the $C^\infty$ Hermitian fibre metric $h\,e^{-f}$ on $E$ is $(\omega,\,\Omega)$-Hermite-Einstein with Einstein factor $c\in\R$. \hfill $\Box$

\vspace{2ex}

In the same vein, one proves the existence of an $(\omega,\,\Omega)$-Hermite-Einstein fibre metric whenever the given vector bundle has rank one, as shown in

\begin{Prop}\label{Prop:H-E_line-bundle}  Let $(X,\,\omega,\,\Omega)$ be a compact K\"ahler manifold equipped with an extra form $\Omega\in C^\infty_{n-m,\,n-m}(X,\,\R)$ satisfying conditions (\ref{eqn:hypotheses_Omega}) for some $m\in\{1,\dots , n\}$. Let $L$ be a holomorphic {\bf line bundle} equipped with a $C^\infty$ Hermitian fibre metric $h$ over $X$.

  Then, there exists a $C^\infty$ function $f:X\longrightarrow\R$, unique up to an additive real constant, such that the rescaled fibre metric $h\,e^{-f}$ is {\bf $(\omega,\,\Omega)$-Hermite-Einstein}, namely there exists a constant $\lambda_{h,\,f}\in\R$ such that \begin{eqnarray}\label{eqn:H-E_prop}i\Theta_{h\,e^{-f}}(L)\wedge\omega^{m-1}\wedge\Omega =  \lambda_{h,\,f}\,dV_\omega.\end{eqnarray}

\end{Prop}

\noindent {\it Proof.} Since $i\Theta_h(L)\wedge\omega^{m-1}\wedge\Omega$ and $dV_\omega$ are $C^\infty$ forms of bidegree $(n,\,n)$ on $X$, there exists a function $\lambda_h:X\longrightarrow\R$ such that \begin{eqnarray*}i\Theta_h(L)\wedge\omega^{m-1}\wedge\Omega = \lambda_h\,dV_\omega\end{eqnarray*} on $X$. Meanwhile, for every function $f:X\longrightarrow\R$, one has $i\Theta_{h\,e^{-f}}(L) = i\Theta_h(L) + i\partial\bar\partial f$, so the rest of the proof follows the same route as that of Proposition \ref{Prop:equiv_w-H-E_H-E}.  \hfill $\Box$

\vspace{2ex}

Again as in the standard case (see e.g. Lemma 2.1.4 in [LT95]), the transformations of the curvature form under the usual operations on vector bundles imply the following

\begin{Lem}\label{Lem:bundle-operations} Let $E,F\longrightarrow(X,\,\omega,\,\Omega)$ be holomorphic vector bundles over a compact K\"ahler manifold equipped with an extra form $\Omega\in C^\infty_{n-m,\,n-m}(X,\,\R)$ satisfying conditions (\ref{eqn:hypotheses_Omega}) for some $m\in\{1,\dots , n\}$. Suppose there exist $(\omega,\,\Omega)$-Hermite-Einstein fibre metrics $h_E$ with Einstein factor $\lambda_E$ on $E$, respectively $h_F$ with Einstein factor $\lambda_F$ on $F$.

  Then, the induced fibre metric on:

\vspace{1ex}

(i)\, $E^\star$ is $(\omega,\,\Omega)$-Hermite-Einstein with Einstein factor $-\lambda_E$;

\vspace{1ex}

(ii)\, $E\otimes F$ is $(\omega,\,\Omega)$-Hermite-Einstein with Einstein factor $\lambda_E + \lambda_F$;

\vspace{1ex}

(iii)\, $\mbox{End}\,(E)$ is $(\omega,\,\Omega)$-Hermite-Einstein with Einstein factor $0$;

\vspace{1ex}

(iv)\, $\Lambda^pE$ is $(\omega,\,\Omega)$-Hermite-Einstein with Einstein factor $p\,\lambda_E$ for every $p\in\{1,\dots , \rk_E\}$;  

\vspace{1ex}

(v)\, $\det E$ is $(\omega,\,\Omega)$-Hermite-Einstein with Einstein factor $\rk_E\,\lambda_E$.

\end{Lem}

\vspace{2ex}

Kobayashi's classical vanishing theorem ([Kob87, Theorem 3.1.9]) takes on the following shape in our setting.

\begin{The}\label{The:vanishing_kobayashi-style} Let $(X,\,\omega,\,\Omega)$ be a compact K\"ahler manifold equipped with an extra form $\Omega\in C^\infty_{n-m,\,n-m}(X,\,\R)$ satisfying conditions (\ref{eqn:hypotheses_Omega}) for some $m\in\{1,\dots , n\}$. Let $E$ be a Hermitian holomorphic vector bundle on $X$ and suppose there exists an {\bf $(\omega,\,\Omega)$-Hermite-Einstein} fibre metric $h$ on $E$ with Einstein factor $\lambda_h$.

\vspace{1ex}

(i)\, If $\lambda_h<0$, then $H^0(X,\,E) = \{0\}$.

\vspace{1ex}

(ii)\, If $\lambda_h=0$, then for every $s\in H^0(X,\,E)$ we have $D_h s=0$ on $X$ (i.e. $s$ is parallel with respect to the Chern connection $D_h$ of $(E,\,h)$).
  
\end{The}

\noindent {\it Proof.} Let $s\in H^0(X,\,E)$. We have: \begin{eqnarray*}P_{\omega,\,\Omega}|s|^2_h = \frac{i\bar\partial\partial\{s,\,s\}_h\wedge\omega^{m-1}\wedge\Omega}{dV_\omega}.\end{eqnarray*} Now, the $h$-compatibility of $D_h$ gives the first inequality below, while the holomorphicity of $s$ yields the second: \begin{eqnarray*}\partial\{s,\,s\}_h = \{D_h's,\,s\}_h + \{s,\,\bar\partial s\}_h = \{D_h's,\,s\}_h.\end{eqnarray*} Taking $\bar\partial$, we further get: \begin{eqnarray*}\bar\partial\partial\{s,\,s\}_h = \{\bar\partial D_h's,\,s\}_h - \{D_h's,\,D_h's\}_h = \{\Theta_h(E)s,\,s\}_h - \{D_h's,\,D_h's\}_h.\end{eqnarray*} Hence: \begin{eqnarray*}i\bar\partial\partial\{s,\,s\}_h\wedge\omega^{m-1}\wedge\Omega & = & \bigg\{\bigg(i\Theta_h(E)\wedge\omega^{m-1}\wedge\Omega\bigg)s,\,s\bigg\}_h - i\,\{D_h's,\,D_h's\}_h\wedge\omega^{m-1}\wedge\Omega \\
  & = & \lambda_h\,\{s,\,s\}_h\,dV_\omega - i\,\{D_h's,\,D_h's\}_h\wedge\omega^{m-1}\wedge\Omega.\end{eqnarray*}

We have thus proved that, for every $s\in H^0(X,\,E)$, we have: \begin{eqnarray}\label{eqn:P_s-squared_formula}P_{\omega,\,\Omega}|s|^2_h = \lambda_h\,|s|_h^2 - \frac{i\,\{D_h's,\,D_h's\}_h\wedge\omega^{m-1}\wedge\Omega}{dV_\omega} \leq  \lambda_h\,|s|_h^2,\end{eqnarray} where the inequality follows by taking $\eta:= D_h's$ in Lemma \ref{Lem:non-negativity_eta_squared} below.

In particular, if $\lambda_h\leq 0$, we get $P_{\omega,\,\Omega}|s|^2_h\leq 0$ at every point of $X$ for every $s\in H^0(X,\,E)$. Since $P_{\omega,\,\Omega}$ is an elliptic second-order differential operator with no zero$^{th}$-order terms and $X$ is compact, the maximum principle implies that \begin{eqnarray*}P_{\omega,\,\Omega}|s|^2_h\equiv 0 \hspace{3ex}\mbox{and}\hspace{3ex} |s|^2_h = Const \hspace{1ex} \mbox{on}\hspace{1ex} X.\end{eqnarray*}

We conclude that:

\vspace{1ex}

$\bullet$ if $\lambda_h< 0$, then $s\equiv 0$ for every $s\in H^0(X,\,E)$;

\vspace{1ex}

$\bullet$ if $\lambda_h= 0$, then $i\,\{D_h's,\,D_h's\}_h\wedge\omega^{m-1}\wedge\Omega \equiv 0$, which amounts to $D_h's = 0$ (hence to $D_hs = 0$ since $\bar\partial s=0$ by holomorphicity of $s$) thanks to Lemma \ref{Lem:non-negativity_eta_squared} applied to $\eta:= D_h's$ for any $s\in H^0(X,\,E)$.

\vspace{1ex}

This is precisely what was claimed.  \hfill $\Box$

\vspace{2ex}

In the above proof, we used the following

\begin{Lem}\label{Lem:non-negativity_eta_squared} Let $(E,\,h)\longrightarrow (X,\,\omega)$ be a Hermitian holomorphic vector bundle with $\rk_\C E = r\geq 1$ over an $n$-dimensional complex Hermitian manifold. Suppose there exist $m\in\{1,\dots , n\}$ and a form $\Omega\in C^\infty_{n-m,\,n-m}(X,\,\C)$ such that $\Omega>0$ (metrically weakly) at every point of $X$.

  Then, for every $\eta\in C^\infty_{1,\,0}(X,\,E)$, the $\R$-valued $C^\infty$ $(n,\,n)$-form $i\,\{\eta,\,\eta\}_h\wedge\omega^{m-1}\wedge\Omega \geq 0$ at every point of $X$. Moreover, for every $x\in X$, we have the equivalence: \begin{eqnarray*}i\,\{\eta,\,\eta\}_h\wedge\omega^{m-1}\wedge\Omega = 0 \hspace{2ex}\mbox{at}\hspace{1ex} x \iff \eta(x)=0.\end{eqnarray*}

\end{Lem}

\noindent {\it Proof.} As in the proof of Lemma \ref{Lem:P-operator_adjoint}, let $\rho$ be the unique $C^\infty$ $(1,\,1)$-form on $X$ that is positive definite and satisfies $\rho_{n-1}:=\rho^{n-1}/(n-1)! = \omega^{m-1}\wedge\Omega$ at every point of $X$. (We use Proposition \ref{Prop:Omega-prod-omega_m-1} and [Mic83] to get existence and uniqueness of this $\rho$.) 

Hence, for every $\eta\in C^\infty_{1,\,0}(X,\,E)$, we get at every point of $X$: \begin{eqnarray}\label{eqn:special-case_norm-eta_10}i\,\{\eta,\,\eta\}_h\wedge\omega^{m-1}\wedge\Omega = i\,\{\eta,\,\eta\}_h\wedge\rho_{n-1} = \Lambda_\rho\bigg(i\,\{\eta,\,\eta\}_h\bigg)\,dV_\rho = |\eta|^2_{\rho,\,h}\,dV_\rho,\end{eqnarray} where $ |\eta|^2_{\rho,\,h}$ is the pointwise squared norm of $\eta$ with respect to $\rho$ and $h$, while  $dV_\rho:=\rho^n/n!$.

This proves the contention once we have proved the last identity in (\ref{eqn:special-case_norm-eta_10}).

To prove this identity, let $U\subset X$ be an open coordinate patch on which $E$ is trivial, let $z_1,\dots , z_n$ be local holomorphic coordinates on $U$ and let $\{e_1,\dots , e_r\}$ be a holomorphic frame of $(E,\,h)$ on $U$. Then, the restriction of $\eta$ to $U$ is of the shape: $\eta_{|U} = \sum\limits_{\alpha=1}^r\eta_\alpha\otimes e_\alpha$, with $\C$-valued $(1,\,0)$-forms $\eta_\alpha$ on $U$. We get at every point of $U$: \begin{eqnarray*}|\eta|^2_{\rho,\,h} = \sum\limits_{\alpha,\,\beta=1}^r \langle\eta_\alpha,\,\eta_\beta\rangle_\rho\,\langle e_\alpha,\,e_\beta\rangle_h \hspace{3ex}\mbox{and}\hspace{3ex} \{\eta,\,\eta\}_h = \sum\limits_{\alpha,\,\beta=1}^r \eta_\alpha\wedge\overline\eta_\beta\,\langle e_\alpha,\,e_\beta\rangle_h.\end{eqnarray*}

Meanwhile, each $\overline\eta_\beta$ is a $\C$-valued $(0,\,1)$-form on $U$, hence it is primitive. Therefore, a well-known formula for the Hodge star operator of any Hermitian metric $\rho$ evaluated on any primitive form yields $\star_\rho\overline\eta_\beta = i\,\overline\eta_\beta\wedge\rho_{n-1}$ for $(0,\,1)$-forms. Hence, for all $\alpha$ and $\beta$, we get: \begin{eqnarray*}\langle\eta_\alpha,\,\eta_\beta\rangle_\rho\,dV_\rho =\eta_\alpha\wedge\star_\rho\overline\eta_\beta = i\,\eta_\alpha\wedge\overline\eta_\beta\wedge\rho_{n-1}.\end{eqnarray*}

Dividing by $dV_\rho$, we get the following well-known formula of independent interest: \begin{eqnarray}\label{eqn:extra-formula_indep}\Lambda_\rho\bigg(i\eta_\alpha\wedge\overline\eta_\beta\bigg) = \langle\eta_\alpha,\,\eta_\beta\rangle_\rho\end{eqnarray} for all $\C$-valued $(1,\,0)$-forms $\eta_\alpha$ and $\eta_\beta$ and any Hermitian metric $\rho$.

  The last identity in (\ref{eqn:special-case_norm-eta_10}) follows by putting together these pieces of information.     \hfill $\Box$

\subsection{Generalised slope stability of torsion-free coherent sheaves}\label{subsection:stability-gen}

We now propose the following generalisations of the classical notions of degree, slope and slope stability for coherent sheaves by replacing again $\omega^{n-m}$ with $\Omega$ in the set-up described at the beginning of this section.


\begin{Def}\label{Def:stability-gen} Let $(X,\,\omega,\,\Omega)$ be a compact K\"ahler manifold equipped with an extra form $\Omega\in C^\infty_{n-m,\,n-m}(X,\,\R)$ satisfying conditions (\ref{eqn:hypotheses_Omega}) for some $m\in\{1,\dots , n\}$.

\vspace{1ex}

(i)\, The {\bf $(\omega,\,\Omega)$-degree} of a holomorphic line bundle $L$ over $X$ is the real number \begin{eqnarray}\label{eqn:degree-L_gen}\deg_{\omega,\,\Omega}(L):=\int\limits_Xi\Theta_h(L)\wedge\omega^{m-1}\wedge\Omega,\end{eqnarray} where $h$ is a Hermitian fibre metric on $L$ and $\Theta_h(L)\in C^\infty_{1,\,1}(X,\,\C)$ is the curvature form of its Chern connection.

\vspace{1ex}

(ii)\, The {\bf $(\omega,\,\Omega)$-degree} of a coherent sheaf ${\cal F}$ of ${\cal O}_X$-modules on $X$ is the real number \begin{eqnarray}\label{eqn:degree-sheaf_gen}\deg_{\omega,\,\Omega}({\cal F}):=\deg_{\omega,\,\Omega}(\det {\cal F}),\end{eqnarray} where $\det {\cal F}$ is the determinant line bundle of ${\cal F}$.

\vspace{1ex}

(iii)\, Let ${\cal F}$ be a coherent, torsion-free non-trivial sheaf of ${\cal O}_X$-modules on $X$. The {\bf $(\omega,\,\Omega)$-slope} of ${\cal F}$ is the real number \begin{eqnarray}\label{eqn:slope_gen}\mu_{\omega,\,\Omega}({\cal F}):=\frac{\deg_{\omega,\,\Omega}({\cal F})}{\rk {\cal F}}.\end{eqnarray}

Furthermore,  ${\cal F}$ is said to be {\bf $(\omega,\,\Omega)$-stable} (resp. {\bf $(\omega,\,\Omega)$-semi-stable}) if for every coherent subsheaf ${\cal G}\subset{\cal F}$ with $0<\rk {\cal G}<\rk {\cal F}$, one has \begin{eqnarray*}\mu_{\omega,\,\Omega}({\cal G})<\mu_{\omega,\,\Omega}({\cal F})  \hspace{5ex} \bigg(\mbox{resp.} \hspace{2ex} \mu_{\omega,\,\Omega}({\cal G})\leq\mu_{\omega,\,\Omega}({\cal F})\bigg).\end{eqnarray*}

\end{Def}

Note that the above definition is correct in the sense that the $(\omega,\,\Omega)$-degree of $L$ is independent of the choice of a Hermitian fibre metric $h$ on $L$. Indeed, for any $h$, the curvature form $i\Theta_h(L)$ belongs to a same Bott-Chern cohomology class $c_1(L)_{BC}\in H^{1,\,1}_{BC}(X,\,\R)$, (the Bott-Chern version of) the first Chern class of $L$. Since $\omega$ is K\"ahler and $\partial\bar\partial\Omega = 0$, we have $\partial\bar\partial(\omega^{m-1}\wedge\Omega) = 0$, so for every $L^1_{loc}$ function $\varphi:X\longrightarrow\R\cup\{-\infty\}$, the compactness of $X$ and the Stokes theorem yield: \begin{eqnarray*}\int\limits_X \bigg(i\Theta_h(L) + i\partial\bar\partial\varphi\bigg)\wedge\omega^{m-1}\wedge\Omega = \int\limits_X i\Theta_h(L)\wedge\omega^{m-1}\wedge\Omega.\end{eqnarray*}

\vspace{2ex}

A quick link with $\S$\ref{subsection:H-E_gen} is given in

\begin{Lem}\label{Lem:Einstein-factor-slope_link} Let $E\longrightarrow (X,\,\omega,\,\Omega)$ be a holomorphic vector bundle over a compact K\"ahler manifold equipped with an extra form $\Omega\in C^\infty_{n-m,\,n-m}(X,\,\R)$ satisfying (\ref{eqn:hypotheses_Omega}). Suppose there exists an {\bf $(\omega,\,\Omega)$-Hermite-Einstein} fibre metric $h$ on $E$.

  Then, the Einstein factor $\lambda_h$ of $h$ is related in the following way to the $(\omega,\,\Omega)$-slope of $E$: \begin{eqnarray}\label{eqn:Einstein-factor-slope_link}\lambda_h = \frac{1}{\mbox{Vol}_\omega(X)}\,\mu_{\omega,\,\Omega}(E),\end{eqnarray} where $\mbox{Vol}_\omega(X):=\int\limits_X\omega^n/n!$ is the $\omega$-volume of $X$.

\end{Lem}
  
\noindent {\it Proof.} Let $r$ be the rank of $E$. From (iii) of Definition \ref
          {Def:stability-gen} we infer the first equality below: \begin{eqnarray*}\mu_{\omega,\,\Omega}(E) & = & \frac{1}{r}\,\int\limits_X i\Theta_h(\det E)\wedge\omega^{m-1}\wedge\Omega = \frac{1}{r}\,\int\limits_X \mbox{Trace}_E\bigg(i\Theta_h(E)\wedge\omega^{m-1}\wedge\Omega\bigg) \\
            & = & \frac{1}{r}\,\int\limits_X \mbox{Trace}_E\bigg(\lambda_h\,dV_\omega\otimes\mbox{Id}_E\bigg) = \frac{1}{r}\,r\,\lambda_h\,\int\limits_X dV_\omega = \lambda_h\,\mbox{Vol}_\omega(X).\end{eqnarray*} This proves the contention. We have used the well-known formula expressing the curvature form of the determinant line bundle $\det E$ as the trace $\mbox{Trace}_E$ in the endomorphism bundle $\mbox{End}\,(E)$ of the curvature form $\Theta_h(E)\in C^\infty_{1,\,1}(X,\,\mbox{End}\,(E))$ of $E$. \hfill $\Box$

\vspace{2ex}

In the setting of Definition \ref{Def:stability-gen}, we consider the function: \begin{eqnarray}\label{eqn:function-degree}d _{\omega,\,\Omega,\,h}({\cal F}):X\longrightarrow\R, \hspace{5ex}  d _{\omega,\,\Omega,\,h}({\cal F}):=\frac{i\Theta_h(\det{\cal F})\wedge\omega^{m-1}\wedge\Omega}{dV_\omega},\end{eqnarray} where $h$ is a Hermitian fibre metric on the determinant line bundle $\det{\cal F}$ of any given torsion-free coherent sheaf ${\cal F}$ on $X$. The following statement generalises Proposition 2.3.1 of [LT95] to our setting.

\begin{Prop}\label{Prop:exact-seq_map_non-decreasing} Let $(X,\,\omega,\,\Omega)$ be a complex K\"ahler manifold equipped with an extra form $\Omega\in C^\infty_{n-m,\,n-m}(X,\,\R)$ satisfying conditions (\ref{eqn:hypotheses_Omega}) for some $m\in\{1,\dots , n\}$. Let \begin{eqnarray}\label{eqn:exact-seq_Hermitian}0\longrightarrow (S,\,h_S) \longrightarrow (E,\,h_E) \longrightarrow (Q,\,h_Q)\longrightarrow 0\end{eqnarray} be a short exact sequence of holomorphic vector bundles on $X$, where $h_E$ is a Hermitian fibre metric on $E$ and $h_S$, resp. $h_Q$, are the induced, resp. quotient, Hermitian fibre metrics on $S$, resp. $Q$.

    If $h_E$ is {\bf $(\omega,\,\Omega)$-Hermite-Einstein} with Einstein factor $\lambda_E$, then: \begin{eqnarray}\label{eqn:exact-seq_map_non-decreasing}\frac{d _{\omega,\,\Omega,\,h_S}(S)}{\rk S}\leq \frac{d _{\omega,\,\Omega,\,h_E}(E)}{\rk E}\leq \frac{d_{\omega,\,\Omega,\,h_Q}(Q)}{\rk Q}.\end{eqnarray}

    Moreover, equality holds in the first inequality if and only if equality holds in the second inequality if and only if the short exact sequence (\ref{eqn:exact-seq_Hermitian}) splits holomorphically. In this case, the Hermitian fibre metrics $h_S$ and $h_Q$ are {\bf $(\omega,\,\Omega)$-Hermite-Einstein} with the same Einstein constants as for $(E,\,h_E)$: \begin{eqnarray*}\lambda_S = \lambda_E = \lambda_Q.\end{eqnarray*}

    In particular, if $X$ is compact, then \begin{eqnarray}\label{eqn:exact-seq_map_non-decreasing-slope}\mu_{\omega,\,\Omega}(S)\leq\mu_{\omega,\,\Omega}(E)\leq\mu_{\omega,\,\Omega}(Q)\end{eqnarray} and equalities occur under the same circumstances as above.

\end{Prop}  
          
\noindent {\it Proof.} From the general theory, we know that the curvature forms of $(S,\,h_S)$ and $(Q,\,h_Q)$ are given in terms of the curvature form $\Theta_{h_E}(E)$ of $(E,\,h_E)$ by the formulae:  \begin{eqnarray}\label{eqn:curvature_S-Q}i\Theta_{h_S}(S) = i\Theta_{h_E}(E)_{|S} + i\beta^\star\wedge\beta \hspace{5ex}\mbox{and}\hspace{5ex} i\Theta_{h_Q}(Q) = i\Theta_{h_E}(E)_{|Q} + i\beta\wedge\beta^\star,\end{eqnarray} where $\beta\in C^\infty_{1,\,0}(X,\,\mbox{Hom}\,(S,\,Q))$ is the second fundamental form of $S$ in $E$ and $\beta^\star\in C^\infty_{0,\,1}(X,\,\mbox{Hom}\,(Q,\,S))$ is its adjoint. Moreover, we always have $\bar\partial\beta^\star = 0$ and the short exact sequence (\ref{eqn:exact-seq_Hermitian}) splits holomorphically if and only if $\beta^\star$ is $\bar\partial$-exact (which further amounts to the extension class $\{\beta^\star\}\in H^{0,\,1}_{\bar\partial}(X,\,\mbox{Hom}\,(Q,\,S))$ vanishing).

Let $s:=\rk S$, $r:=\rk E$ and $q:=\rk Q$. We get: \begin{eqnarray}\label{eqn:d_E_computation}\nonumber d _{\omega,\,\Omega,\,h_E}(E) & = & \frac{i\Theta_{h_E}(\det E)\wedge\omega^{m-1}\wedge\Omega}{dV_\omega} = \mbox{Trace}_E\bigg(\frac{i\Theta_{h_E}(E)\wedge\omega^{m-1}\wedge\Omega}{dV_\omega}\bigg) \\
  & = & \mbox{Trace}_E\bigg(\lambda_E\,\mbox{Id}_E\bigg) = r\,\lambda_E,\end{eqnarray} where the last but one equality followed from the $(\omega,\,\Omega)$-Hermite-Einstein hypothesis on $(E,\,h_E)$.

Similarly, we get (using also (\ref{eqn:curvature_S-Q}) to get the first equality): \begin{eqnarray}\label{eqn:d_S_computation}d _{\omega,\,\Omega,\,h_S}(S) & = & \mbox{Trace}_S\bigg(\frac{i\Theta_{h_E}(E)_{|S}\wedge\omega^{m-1}\wedge\Omega}{dV_\omega}\bigg) + \mbox{Trace}_S\bigg(\frac{i\beta^\star\wedge\beta\wedge\omega^{m-1}\wedge\Omega}{dV_\omega}\bigg) \\
\nonumber & = & \mbox{Trace}_S\bigg(\lambda_E\,\mbox{Id}_S\bigg) + \mbox{Trace}_S\bigg(\frac{i\beta^\star\wedge\beta\wedge\omega^{m-1}\wedge\Omega}{dV_\omega}\bigg) = s\lambda_E + \mbox{Trace}_S\bigg(\frac{i\beta^\star\wedge\beta\wedge\omega^{m-1}\wedge\Omega}{dV_\omega}\bigg)\end{eqnarray} and \begin{eqnarray}\label{eqn:d_Q_computation}d _{\omega,\,\Omega,\,h_Q}(Q) = q\lambda_E + \mbox{Trace}_Q\bigg(\frac{i\beta\wedge\beta^\star\wedge\omega^{m-1}\wedge\Omega}{dV_\omega}\bigg)\end{eqnarray}

Taking (\ref{eqn:d_E_computation}), (\ref{eqn:d_S_computation}) and (\ref{eqn:d_Q_computation}) into account, we see that proving (\ref{eqn:exact-seq_Hermitian}) is equivalent to proving: \begin{eqnarray*}\mbox{Trace}_S\bigg(\frac{i\beta^\star\wedge\beta\wedge\omega^{m-1}\wedge\Omega}{dV_\omega}\bigg)\leq 0\leq\mbox{Trace}_Q\bigg(\frac{i\beta\wedge\beta^\star\wedge\omega^{m-1}\wedge\Omega}{dV_\omega}\bigg),\end{eqnarray*} a fact that is guaranteed by Lemma \ref{Lem:non-negativity_beta-end_squared} below.

This proves the first statement. The other statements follow at once from this one. \hfill $\Box$

\vspace{2ex}

In the above proof, we used the following analogue of Lemma \ref{Lem:non-negativity_eta_squared}.

\begin{Lem}\label{Lem:non-negativity_beta-end_squared} The setting is the same as in Proposition \ref{Prop:exact-seq_map_non-decreasing} except that we do not make any $(\omega,\,\Omega)$-Hermite-Einstein assumption on $(E,\,h_E)$.

  Then, the forms $i\beta\wedge\beta^\star\in C^\infty_{1,\,1}(X,\,\mbox{Hom}\,(Q,\,Q))$ and $i\beta^\star\wedge\beta\in C^\infty_{1,\,1}(X,\,\mbox{Hom}\,(S,\,S))$, constructed from the second fundamental form $\beta\in C^\infty_{1,\,0}(X,\,\mbox{Hom}\,(S,\,Q))$ of $S$ in $E$ and its adjoint $\beta^\star\in C^\infty_{0,\,1}(X,\,\mbox{Hom}\,(Q,\,S))$, satisfy the following inequalities: \begin{eqnarray*}\mbox{Trace}_Q\bigg(\frac{i\beta\wedge\beta^\star\wedge\omega^{m-1}\wedge\Omega}{dV_\omega}\bigg) \geq 0 \hspace{3ex}\mbox{and}\hspace{3ex} \mbox{Trace}_S\bigg(\frac{i\beta^\star\wedge\beta\wedge\omega^{m-1}\wedge\Omega}{dV_\omega}\bigg)\leq 0\end{eqnarray*} at every point of $X$.

    Moreover, for every $x\in X$, either of the above inequalities is an equality at $x$ if and only if $\beta(x) = 0$.

\end{Lem}

\noindent {\it Proof.} On an open coordinate patch $U\subset X$ on which $E$ is trivial, let $z_1,\dots , z_n$ be local holomorphic coordinates on $X$ and let $\{e_1,\dots , e_r\}$ be a frame for $E$ such that $\{e_1,\dots , e_s\}$ is a frame for $S$ and $\{e_{s+1},\dots , e_r\}$ is a frame for $Q$. (Recall that $s+q = r$.) Then, on $U$ we have: \begin{eqnarray*}\beta = \sum\limits_{j=1}^n\sum\limits_{\lambda=1}^s\sum\limits_{\mu=s+1}^r\beta_j^{\lambda\mu}\,dz_j\otimes e_\lambda^\star\otimes e_\mu, \hspace{3ex}\mbox{hence}\hspace{3ex} \beta^\star = \, ^t\overline\beta= \sum\limits_{j=1}^n\sum\limits_{\lambda=1}^s\sum\limits_{\mu=s+1}^r\overline{\beta_j^{\lambda\mu}}\,d\bar{z}_j\otimes e_\mu^\star\otimes e_\lambda.\end{eqnarray*} From this, we get the following erqualities on $U$: \begin{eqnarray*}i\beta\wedge\beta^\star & = & \sum\limits_{j,\,k=1}^n\sum\limits_{\lambda=1}^s\sum\limits_{\mu,\,\nu=s+1}^r \beta_j^{\lambda\nu}\,\overline{\beta_k^{\lambda\mu}}\,idz_j\wedge d\bar{z}_k\otimes e_\nu\otimes e_\mu^\star \\
  \mbox{Trace}_Q(i\beta\wedge\beta^\star)  & = & \sum\limits_{j,\,k=1}^n\sum\limits_{\lambda=1}^s\sum\limits_{\mu=s+1}^r \beta_j^{\lambda\mu}\,\overline{\beta_k^{\lambda\mu}}\,idz_j\wedge d\bar{z}_k.\end{eqnarray*}

On the other hand, as in the proof of Lemma \ref{Lem:P-operator_adjoint}, let $\rho$ be the unique $C^\infty$ $(1,\,1)$-form on $X$ that is positive definite and satisfies $\rho_{n-1}: = \rho^{n-1}/(n-1)! = \omega^{m-1}\wedge\Omega$ at every point of $X$.

Now, fix an arbitrary point $x\in U$ and choose the $z_j$'s on $U$ centred at $x$ such that $\rho(x) = \sum_j idz_j\wedge d\bar{z}_j$. Then, at $x$ we get: \begin{eqnarray*}\mbox{Trace}_Q\bigg(\frac{i\beta\wedge\beta^\star\wedge\omega^{m-1}\wedge\Omega}{dV_\rho}\bigg) = \bigg(\sum\limits_{j=1}^n\sum\limits_{\lambda=1}^s\sum\limits_{\mu=s+1}^r|\beta_j^{\lambda\mu}|^2\bigg)\,dV_n\geq 0, \end{eqnarray*} where $dV_n:=idz_1\wedge d\bar{z}_1\wedge\dots\wedge idz_n\wedge d\bar{z}_n$. Moreover, we see that equality occurs in the last inequality if and only if $\beta_j^{\lambda\mu} = 0$ for all $j,\lambda,\mu$, a fact that amounts to $\beta(x) = 0$.

Similarly, we get \begin{eqnarray*}\mbox{Trace}_S\bigg(\frac{i\beta^\star\wedge\beta\wedge\omega^{m-1}\wedge\Omega}{dV_\rho}\bigg) = -\mbox{Trace}_Q\bigg(\frac{i\beta\wedge\beta^\star\wedge\omega^{m-1}\wedge\Omega}{dV_\rho}\bigg)\leq 0,\end{eqnarray*} with equality if and only if $\beta(x) = 0$.

Since $dV_\rho$ and $dV_\omega$ differ by a factor that is positive at every point of $X$, the contention follows. \hfill $\Box$

\vspace{2ex}

We are now ready to give the main result of this paper. It links the two notions presented above. It is the exact analogue, in our setting, of the standard result established by Kobayashi and L\"ubke [Kob82, Lub83] and stated as Theorem 2.3.2 in [LT95].

\begin{The}\label{The:stability-H-E_implication_gen}  Let $(X,\,\omega,\,\Omega)$ be a compact K\"ahler manifold equipped with an extra form $\Omega\in C^\infty_{n-m,\,n-m}(X,\,\R)$ satisfying conditions (\ref{eqn:hypotheses_Omega}) for some $m\in\{1,\dots , n\}$. Let $E$ be a holomorphic vector bundle of rank $r\geq 1$ over $X$.

  Suppose there exists an {\bf $(\omega,\,\Omega)$-Hermite-Einstein} $C^\infty$ fibre metric $h$ on $E$. Then $E$ is {\bf $(\omega,\,\Omega)$-semi-stable}. If, moreover, $E$ is not $(\omega,\,\Omega)$-stable, it is a holomorphic, $h$-orthogonal direct sum \begin{eqnarray*}E = \bigoplus\limits_i E_i\end{eqnarray*} of $(\omega,\,\Omega)$-Hermite-Einstein holomorphic subbundles $E_i\subset E$, the fibre metric $h_i$ induced by $h$ in each $E_i$ being {\bf $(\omega,\,\Omega)$-Hermite-Einstein} with the same Einstein factor as $(E,\,h)$.

\end{The}

\noindent {\it Proof.} Let ${\cal F}\subset{\cal O}(E)$ be a coherent subsheaf of some rank $s$ with $0<s<r$. The inclusion $\varphi:{\cal F}\hookrightarrow{\cal O}(E)$ induces a holomorphic morphism $\det\varphi:\det{\cal F}\hookrightarrow\bigg(\Lambda^s{\cal O}(E)\bigg)^{\star\star}\simeq\Lambda^s{\cal O}(E)$ which is injective because ${\cal F}$ is locally free on a Zariski open subset of $X$ (actually, outside an analytic subset of codimension $\geq 2$). By tensoring with $\det{\cal F}^\star$ (which is locally free of rank $1$), we get a non-trivial global holomorphic section $\sigma\in H^0(X,\,\Lambda^s E\otimes\det{\cal F}^\star)$: \begin{eqnarray}\label{eqn:global-section}\sigma:\det{\cal F}\otimes\det{\cal F}^\star\simeq{\cal O}_X \hookrightarrow \Lambda^s E\otimes\det{\cal F}^\star.\end{eqnarray}   

Let $\lambda_h$ be the Einstein factor of $h$, so \begin{eqnarray*}i\Theta_h(E)\wedge\omega^{m-1}\wedge\Omega =  \lambda_h\,dV_\omega\otimes\mbox{Id}_E.\end{eqnarray*} Moreover, $\lambda_h$ is given in terms of the $(\omega,\,\Omega)$-slope $\mu_{\omega,\,\Omega}(E)$ of $E$ by formula (\ref{eqn:Einstein-factor-slope_link}). Meanwhile, $\det{\cal F}$ is a holomorphic line bundle on a compact $X$, so it carries an {\bf $(\omega,\,\Omega)$-Hermite-Einstein} fibre metric $h_{\det{\cal F}}$ by Proposition \ref{Prop:H-E_line-bundle}. By the same formula (\ref{eqn:Einstein-factor-slope_link}), its Einstein factor $\lambda(\det{\cal F})\in\R$ satisfies: \begin{eqnarray*}\lambda(\det{\cal F}) = \frac{1}{\mbox{Vol}_\omega(X)}\,\mu_{\omega,\,\Omega}(\det{\cal F}) = \frac{1}{\mbox{Vol}_\omega(X)}\,\deg_{\omega,\,\Omega}({\cal F}) = \frac{1}{\mbox{Vol}_\omega(X)}\,\mu_{\omega,\,\Omega}({\cal F}).\end{eqnarray*}

Hence, the metric induced on $\Lambda^s E\otimes\det{\cal F}^\star$ by $h$ and $h_{\det{\cal F}}$ is, thanks to Lemma \ref{Lem:bundle-operations}, {\bf $(\omega,\,\Omega)$-Hermite-Einstein} with Einstein factor \begin{eqnarray*}s\,\lambda_h - \lambda(\det{\cal F}) = \frac{s}{\mbox{Vol}_\omega(X)}\,\bigg(\mu_{\omega,\,\Omega}(E) - \mu_{\omega,\,\Omega}({\cal F})\bigg).\end{eqnarray*}

On the other hand, since the holomorphic $(\omega,\,\Omega)$-Hermite-Einstein vector bundle $\Lambda^s E\otimes\det{\cal F}^\star$ has a non-trivial global holomorphic section (see (\ref{eqn:global-section})), our version of Kobayashi's vanishing theorem (cf. Theorem \ref{The:vanishing_kobayashi-style}) implies that $s\,\lambda_h - \lambda(\det{\cal F})\geq 0$, a fact that amounts to \begin{eqnarray*}\mu_{\omega,\,\Omega}({\cal F})\leq\mu_{\omega,\,\Omega}(E)\end{eqnarray*} thanks to what was concluded above. This proves that $E$ is $(\omega,\,\Omega)$-semi-stable.

\vspace{2ex}

Now, suppose that $E$ is not $(\omega,\,\Omega)$-stable. Then, there exists a coherent subsheaf ${\cal F}\subset{\cal O}(E)$ of some rank $s$ with $0<s<r$ such that $\mu_{\omega,\,\Omega}({\cal F}) = \mu_{\omega,\,\Omega}(E)$, a fact that amounts to $s\,\lambda_h - \lambda(\det{\cal F}) = 0$.

We may assume that ${\cal F}$ is reflexive (or, equivalently, torsion-free and normal). Thanks to our version of Kobayashi's vanishing theorem (cf. Theorem \ref{The:vanishing_kobayashi-style}), the section $\sigma\in H^0(X,\,\Lambda^s E\otimes\det{\cal F}^\star)$ must be parallel, so the holomorphic line bundle $\det{\cal F}$ must be a parallel subbundle of $\Lambda^s E$. This implies that the {\bf $(\omega,\,\Omega)$-Hermite-Einstein} fibre metric $h_{\det{\cal F}}$ we found on $\det{\cal F}$ equals (up to a positive multiplicative constant) the metric induced from the one on $\Lambda^s E$ by the inclusion $\det\varphi:\det{\cal F}\hookrightarrow\Lambda^s E$.

Now, let $S\subset X$ be the singular locus of ${\cal F}$. Then:

\vspace{1ex}

$\bullet$ ${\cal F}_{|Y}$ is locally free, where $Y:=X\setminus S$;

\vspace{1ex}

$\bullet$ $S$ is an analytic subset of $X$ and, since ${\cal F}$ is reflexive, $\mbox{codim}_XS\geq 3$. 

\vspace{1ex}

Meanwhile, the assumption $\mu_{\omega,\,\Omega}({\cal F}) = \mu_{\omega,\,\Omega}(E)$ is equivalent to \begin{eqnarray*}\frac{1}{s}\,\int\limits_Xi\Theta_{h_{\cal F}}(\det {\cal F})\wedge\omega^{m-1}\wedge\Omega = \frac{1}{r}\,\int\limits_Xi\Theta_h(\det E)\wedge\omega^{m-1}\wedge\Omega\end{eqnarray*} and, given the above discussion, this implies the equality of functions: \begin{eqnarray*}\frac{d _{\omega,\,\Omega,\,h_{\cal F}}({\cal F}_{|Y})}{s} = \frac{d _{\omega,\,\Omega,\,h}(E_{|Y})}{r} \hspace{5ex} \mbox{on}\hspace{1ex} Y,\end{eqnarray*} where $h_{\cal F}$ is the metric induced on ${\cal F}$ by $\varphi$ and $h$.

  By Proposition \ref{Prop:exact-seq_map_non-decreasing}, this equality of functions is equivalent to the short exact sequence of holomorphic vector bundles on $Y$ \begin{eqnarray*}0\longrightarrow{\cal F}_{|Y}\longrightarrow E_{|Y}\longrightarrow E_{|Y}/{\cal F}_{|Y} \longrightarrow 0\end{eqnarray*} splitting holomorphically. Moreover, this also implies that ${\cal F}_{|Y}$ and $E_{|Y}/{\cal F}_{|Y}$ are $(\omega,\,\Omega)$-Hermite-Einstein with the same Einstein factor as $(E,\,h)$.

    Since the sheaf $\mbox{Hom} ({\cal O}(E),\,{\cal F})$ is normal and $\mbox{codim}_XS\geq 3$, the splitting morphism ${\cal O}(E)_{|Y}\longrightarrow{\cal F}_{|Y}$ extends to a morphism ${\cal O}(E)\longrightarrow{\cal F}$. Thus, the exact sequence of sheaves \begin{eqnarray*}0\longrightarrow{\cal F}\longrightarrow {\cal O}(E)\longrightarrow {\cal G}:={\cal O}(E)/{\cal F} \longrightarrow 0\end{eqnarray*} splits on $X$. In particular, ${\cal F}$ and ${\cal G}$ are locally free, hence they define holomorphic subbundles $F$ and $G$ of $E$ over $X$ and we have a holomorphic direct-sum decomposition: \begin{eqnarray*}E = F\oplus G.\end{eqnarray*} Moreover, the subbundles $F$ and $G$ of $E$ are $(\omega,\,\Omega)$-Hermite-Einstein with the same Einstein factor as $(E,\,h)$.

An induction on the rank of $E$ completes the proof. \hfill $\Box$

\section{Appendix}\label{section:appendix} The notions of weak and strong semi-positive $(p,\,p)$-forms and currents on a complex manifold $X$ and on a given $\C$-vector space $V$ (that is usually taken to be the holomorphic tangent space $T_x^{1,\,0}X$ to a complex manifold $X$ at a point $x\in X$ in geometric situations) are standard. They go back to Lelong (1957) and are given a detailed exposition in $\S1$ of chapter III of Demailly's book [Dem97]. However, we call ``semi-positive'' what is termed ``positive'' in [Dem97] and we reinforce these notions to their strictly positive counterparts which do not seem to have been formalised in the literature so far. Hence this appendix.

\subsection{(Semi-)positivity for forms}\label{subsection:semi_pos_forms} We first deal with differential forms. Currents will be treated in the next subsection.  

\subsubsection{Classical notions for forms: semi-positivity}\label{subsubsection:classical}

For the reader's convenience, we start by recalling the classical Lelong notions of (weak and strong) semi-positivity (cf. [Dem97, III, $\S1.1$]). A $(p,\,p)$-form $\Omega\in\Lambda^{p,\,p}V^\star$ on an $n$-dimensional $\C$-vector space $V$ is said to be:

\vspace{1ex}

(i)\, {\bf weakly semi-positive} if for all $(1,\,0)$-forms $\alpha_1,\dots , \alpha_{n-p}\in V^\star$ on $V$, the following $(n,\,n)$-form is {\bf non-negative}: \begin{eqnarray*}\Omega\wedge i\alpha_1\wedge\bar\alpha_1\wedge\dots\wedge  i\alpha_{n-p}\wedge\bar\alpha_{n-p}\geq 0;\end{eqnarray*} We write in this case: $\Omega\geq 0$ {\bf (weakly)}.

\vspace{1ex}

(ii)\, {\bf strongly semi-positive} if $\Omega$ is a {\it finite} linear combination with {\bf non-negative} coefficients of {\bf decomposable} $(p,\,p)$-forms, namely if there exist $(1,\,0)$-forms $\beta_{j,\,1},\dots , \beta_{j,\,p}\in V^\star$ on $V$, with $j$ varying in a finite set of indices, and corresponding reals $\lambda_j\geq 0$ such that \begin{eqnarray*}\Omega = \sum\limits_j\lambda_j\,i\beta_{j,\,1}\wedge\bar\beta_{j,\,1}\wedge\dots\wedge i\beta_{j,\,p}\wedge\bar\beta_{j,\,p}.\end{eqnarray*} We write in this case: $\Omega\geq 0$ {\bf (strongly)}.

These definitions naturally lead to their manifold counterparts. A $(p,\,p)$-form $\Omega$ on an $n$-dimensional complex manifold $X$ is said to be {\bf weakly semi-positive} (respectively {\bf strongly semi-positive}) if for every point $x\in X$, the $(p,\,p)$-form $\Omega(x)\in\Lambda^{p,\,p} T_x^\star X$ is {\it weakly semi-positive} (respectively {\it strongly semi-positive}) in the above sense with $V=T_x^{1,\,0}X$.

\vspace{1ex}

Moreover, these definitions imply that, for every $p$, the cone of {\it weakly semi-positive} $(p,\,p)$-forms is dual to the cone of {\it strongly semi-positive} $(n-p,\,n-p)$-forms via the bilinear pairing: \begin{eqnarray}\label{eqn:duality-pos-pairing}\Lambda^{p,\,p}V^\star\times\Lambda^{n-p,\,n-p}V^\star\longrightarrow\C, \hspace{5ex} (u,\,v)\longmapsto\frac{u\wedge v}{d\mu},\end{eqnarray} where $d\mu>0$ is any fixed positive $(n,\,n)$-form on $V$. In other words, a given $u\in\Lambda^{p,\,p}V^\star$ is {\it weakly semi-positive} if and only if $u\wedge v\geq 0$ for all {\it strongly semi-positive} forms $v\in\Lambda^{n-p,\,n-p}V^\star$. Similarly, a given $u\in\Lambda^{p,\,p}V^\star$ is {\it strongly semi-positive} if and only if $u\wedge v\geq 0$ for all {\it weakly semi-positive} forms $v\in\Lambda^{n-p,\,n-p}V^\star$.

This last statement needs a little more context. The previous one follows at once from the definition of weakly semi-positive forms and is equivalent to the fact that the cone of {\it weakly semi-positive} $(p,\,p)$-forms equals the dual of the cone of {\it strongly semi-positive} $(n-p,\,n-p)$-forms. Taking duals, we infer that the {\it dual} of the cone of {\it weakly semi-positive} $(p,\,p)$-forms equals the {\it bidual} of the cone of {\it strongly semi-positive} $(n-p,\,n-p)$-forms. On the other hand, it is standard that the {\it bidual} of any convex cone $\Gamma$ equals its closure $\overline\Gamma$. So, it remains to prove that the cone of {\it strongly semi-positive} $(n-p,\,n-p)$-forms is {\it closed}. The analogous statement for the cone of {\it weakly semi-positive} $(p,\,p)$-forms is obvious. What might a priori prevent the former cone from being closed is the following possibility.\footnote{The author is grateful to B. Berndtsson for pointing out this aspect to him.} If $(\Omega_\nu)_{\nu\geq 1}$ is a convergent sequence of {\it strongly semi-positive} $(n-p,\,n-p)$-forms, the number of terms (each being a non-negative multiple of a decomposable form) in the finite sum defining $\Omega_\nu$ may increase and tend to $\infty$ as $\nu\to\infty$. A priori, one needs to justify that the limiting form is still a {\it finite} linear combination with non-negative coefficients of decomposable forms. We give a proof of the closedness of the cone of {\it strongly semi-positive} $(p,\,p)$-forms, for any $p$, in Lemma \ref{Lem:s-semi-pos-cone_closed} below.

\vspace{2ex}

Other standard facts (see [Dem97, III, $\S1.1$]) are:

\vspace{1ex}

(a)\, every {\it strongly semi-positive} $(p,\,p)$-form is {\it weakly semi-positive};

\vspace{1ex}

(b)\, every {\it weakly semi-positive} $(p,\,p)$-form $\Omega$ is {\it real} (in the sense that $\overline\Omega = \Omega$);

\vspace{1ex}

(c)\, for $p=1$ and $p=n-1$, the notions of {\it strongly semi-positive} $(p,\,p)$-forms and {\it weakly semi-positive} $(p,\,p)$-forms coincide; a $(1,\,1)$-form $u = \sum_{1\leq j,\,k\leq n} u_{j\bar{k}}\,idz_j\wedge d\bar{z}_k$ is {\it semi-positive} if and only if its coefficient matrix $(u_{j\bar{k}})_{1\leq j,\,k\leq n}$ is a positive semi-definite Hermitian matrix;

\vspace{1ex}

(d)\, for every $p\in\{2,\dots , n-2\}$, there exist {\it weakly semi-positive} $(p,\,p)$-forms that are {\it not strongly semi-positive};

\vspace{1ex}

(e)\, a form $\Omega\in\Lambda^{p,\,p}V^\star$ is {\it weakly semi-positive} if and only if its restriction $\Omega_{|S}$ to every $p$-dimensional $\C$-vector subspace $S\subset V$ satisfies $\Omega_{|S}\geq 0$ as a top-degree form on $S$;

\vspace{1ex}

(f)\, whenever $u_1,\dots , u_s$ are {\it strongly semi-positive} forms, the exterior product $u_1\wedge\dots \wedge u_s$ is {\it strongly semi-positive}; whenever all but one of $u_1,\dots , u_s$ are {\it strongly semi-positive} and the remaining one is {\it weakly semi-positive}, the exterior product $u_1\wedge\dots \wedge u_s$ is {\it weakly semi-positive};

\vspace{1ex}

(g)\, for every $(p,\,0)$-form $\Gamma\in\Lambda^{p,\,0}V^\star$, the $(p,\,p)$-form $i^{p^2}\,\Gamma\wedge\overline\Gamma$ is {\it weakly semi-positive}; when $2\leq p\leq n-2$, $i^{p^2}\,\Gamma\wedge\overline\Gamma$ is {\it strongly semi-positive} if and only if there exist $(1,\,0)$-forms $\beta_1,\dots , \beta_p\in V^\star$ such that $\beta = \beta_1\wedge\dots \wedge\beta_p$.

\vspace{2ex}

We end this discussion of standard material with the following

\begin{Obs}\label{Obs:m-m_diagonal_equiv} Let $V$ be an $n$-dimensional $\C$-vector space and let $p\in\{1,\dots , n\}$. Let $\{\alpha_1,\dots , \alpha_n\}$ be a $\C$-basis of $V^\star$. For every multi-index $I=(1\leq i_1<\dots < i_p\leq n)$ of length $|I|=p$, set $\alpha_I:=\alpha_{i_1}\wedge\dots\wedge\alpha_{i_p}\in\Lambda^{p,\,0}V^\star$.

  Then, for any {\bf diagonal} $(p,\,p)$-form, namely any form of the shape \begin{eqnarray}\label{eqn:diagonal-shape}\Omega = \sum\limits_{|I|=p}\lambda_I\,i^{p^2}\,\alpha_I\wedge\bar\alpha_I = \sum\limits_{|I|=p}\lambda_I\,i\alpha_{i_1}\wedge\bar\alpha_{i_1}\wedge\dots\wedge i\alpha_{i_p}\wedge\bar\alpha_{i_p}\in\Lambda^{p,\,p}V^\star, \hspace{5ex} \lambda_I\in\C,\end{eqnarray} the following equivalences hold: \begin{eqnarray*}\Omega\geq 0 \hspace{2ex} \mbox{(strongly)} \iff \Omega\geq 0 \hspace{2ex} \mbox{(weakly)} \iff \lambda_I\geq 0 \hspace{5ex} \mbox{for every}\hspace{1ex} I  \hspace{1ex} \mbox{with}\hspace{1ex} |I|=p.\end{eqnarray*}

\end{Obs}  

\noindent {\it Proof.} The implication ``$\Omega\geq 0$ (strongly) $\implies$ $\Omega\geq 0$ (weakly)'' is classical (and holds for every $(p,\,p)$-form).

Now, suppose that $\Omega\geq 0$ (weakly). Then, for every $I$ with $|I|=p$, denoting by $C_I:=\{1,\dots , n\}\setminus I$ the multi-index (of length $n-p$) complementary to $I$, the $(n-p,\,n-p)$-form $i^{(n-p)^2}\,\alpha_{C_I}\wedge\bar\alpha_{C_I}$ is decomposable. Therefore, we get the first inequality below: \begin{eqnarray*}0\leq\Omega\wedge i^{(n-p)^2}\,\alpha_{C_I}\wedge\bar\alpha_{C_I} = \lambda_I\, dV_n,\end{eqnarray*} where $dV_n:=i\alpha_1\wedge\bar\alpha_1\wedge\dots\wedge i\alpha_n\wedge\bar\alpha_n>0$ is the volume form induced on $V$ by the given basis of $V^\star$. We infer that $\lambda_I\geq 0$. Since $I$ was arbitrary, the second direct implication in the statement follows.

Finally, if $\lambda_I\geq 0$ for every $I$ with $|I|=p$, then $\Omega$ is a finite linear combination with non-negative coefficients of decomposable $(p,\,p)$-forms. Thus, $\Omega\geq 0$ (strongly) and we are done. \hfill $\Box$

\vspace{2ex}

Since the notions of weak and strong semi-positivity are independent of the choice of $\C$-basis of $V^\star$, the equivalences in Observation \ref{Obs:m-m_diagonal_equiv} also apply to {\bf diagonalisable} $(p,\,p)$-forms (i.e. forms $\Omega$ that have the shape (\ref{eqn:diagonal-shape}) with respect to some $\C$-basis of $V^\star$). In particular, we recover the standard result according to which the notions of weak and strong semi-positivity are equivalent in bidegrees $(1,\,1)$ and $(n-1,\,n-1)$. Indeed, every form is diagonalisable in these bidegrees.   

An immediate consequence of Observation \ref{Obs:m-m_diagonal_equiv} is that the notion of $m$-semi-positivity for $(1,\,1)$-forms considered in [DP25a] and [DP25b] (that had been introduced by Dieu and subsequently studied by various authors, including Harvey and Lawson, Verbitsky, Dinew) can be defined in an equivalent way by means of {\bf strong} semi-positivity or {\bf weak} semi-positivity for $(m,\,m)$-forms.

Specifically, if $(X,\,\omega)$ is a complex Hermitian manifold with $\mbox{dim}_\C X = n$, $m\in\{1,\dots , n\}$ and $\alpha$ is a real (for example $C^\infty$ or just continuous) $(1,\,1)$-form on $X$, the next corollary gives the following equivalences at every point $x\in X$: \begin{eqnarray*}\alpha\wedge\omega^{m-1}\geq 0 \hspace{2ex} \mbox{(strongly)}  \iff \alpha\wedge\omega^{m-1}\geq 0 \hspace{2ex} \mbox{(weakly)} \iff \lambda_1 + \dots + \lambda_m\geq 0,\end{eqnarray*} where $\lambda_1\leq\dots\leq\lambda_m\leq\dots\leq\lambda_n$ are the eigenvalues of $\alpha$ with respect to $\omega$ at $x\in X$. The  $(1,\,1)$-form $\alpha$ is said to be {\bf $m$-semi-positive with respect to $\omega$} if any of the above three equivalent conditions is satisfied at every point $x\in X$.

\begin{Cor}\label{Obs:m-pos} Let $V$ be an $n$-dimensional $\C$-vector space and let $m\in\{1,\dots , n\}$. Fix a positive definite $(1,\,1)$-form $\omega\in\Lambda^{1,\,1}V^\star$ on $V$.

  Then, for every $(1,\,1)$-form $\alpha\in\Lambda^{1,\,1}V^\star$, the following equivalences hold: \begin{eqnarray*}\alpha\wedge\omega^{m-1}\geq 0 \hspace{2ex} \mbox{(strongly)} \iff \alpha\wedge\omega^{m-1}\geq 0 \hspace{2ex} \mbox{(weakly)} \iff \lambda_1 + \dots + \lambda_m\geq 0,\end{eqnarray*} where $\lambda_1\leq\dots\leq\lambda_m\leq\dots\leq\lambda_n$ are the eigenvalues of $\alpha$ with respect to $\omega$.

\end{Cor}

\noindent {\it Proof.} It suffices to choose a $\C$-basis $\{\alpha_1,\dots , \alpha_n\}$ of $V^\star$ with respect to which $\alpha$ and $\omega$ are simultaneously diagonal, namely \begin{eqnarray*}\label{eqn:coordinates_diagonalisation}\omega = \sum\limits_{j=1}^n i\alpha_j\wedge d\bar\alpha_j  \hspace{5ex}\mbox{and}\hspace{5ex}  \alpha = \sum\limits_{j=1}^n \lambda_j\,i\alpha_j\wedge d\bar\alpha_j,\end{eqnarray*} to notice that \begin{eqnarray*}\alpha\wedge\omega_{m-1} = \sum\limits_{|J|=m}\bigg(\sum\limits_{j\in J}\lambda_j\bigg)\,i^{m^2} \alpha_J\wedge\bar\alpha_J,\end{eqnarray*} and to apply Observation \ref{Obs:m-m_diagonal_equiv}. \hfill $\Box$

\subsubsection{Other cones of forms}\label{subsubsection:other-cones}

For context, we now include a discussion of related cones.

Let $V$ be a $\C$-vector space with $\mbox{dim}_\C V =n$. Let $p\in\{1,\dots , n\}$ and set $N:={n \choose p}$. Thus, $\mbox{dim}_\C\Lambda^{p,\,p}V^\star = N^2$. Inside the $\C$-vector space $\Lambda^{p,\,p}V^\star$ sits the $\R$-vector space $(\Lambda^{p,\,p}V^\star)_\R$ of real $(p,\,p)$-forms on $V$: \begin{eqnarray*}(\Lambda^{p,\,p}V^\star)_\R = \bigg\{u\in\Lambda^{p,\,p}V^\star\,\mid\, \bar{u} = u\bigg\}\subset\Lambda^{p,\,p}V^\star.\end{eqnarray*} One has $\mbox{dim}_\R(\Lambda^{p,\,p}V^\star)_\R = N^2$.

\vspace{1ex}

We define three convex cones in $(\Lambda^{p,\,p}V^\star)_\R$ and will subsequently discuss their mutual relations.

\begin{Def}\label{Def:cone_all-decomposable} Consider the following set: \begin{eqnarray}\label{eqn:cone_all-decomposable}{\cal C}:=\bigg\{\Xi = \sum\limits_j \lambda_j\,i\beta_{j,\,1}\wedge\bar\beta_{j,\,1}\wedge\dots\wedge i\beta_{j,\,p}\wedge\bar\beta_{j,\,p} \,\mid\,\lambda_j\geq 0 \hspace{1ex}\mbox{and}\hspace{1ex} \beta_{j,\,1},\dots , \beta_{j,\,p}\in V^\star\bigg\}\subset(\Lambda^{p,\,p}V^\star)_\R,\end{eqnarray} where the sums over $j$ are finite.

\end{Def}  

This is the convex cone generated by {\it all} decomposable $(p,\,p)$-forms, the usual cone of {\it strongly semi-positive} $(p,\,p)$-forms on $V$. It is a {\it closed} convex cone in the finite-dimensional real vector space $(\Lambda^{p,\,p}V^\star)_\R$ and it is {\it intrinsic}, namely independent of the choice of any basis in $V^\star$ or $\Lambda^{p,\,p}V^\star$. Its closedness, long taken for granted in the literature, was also noticed and proved in [Xia26, Lemma A.3.]. The proof we give below seems different.

\begin{Lem}\label{Lem:s-semi-pos-cone_closed} The cone ${\cal C}$ is {\bf closed} in $(\Lambda^{p,\,p}V^\star)_\R$.

\end{Lem}   

\noindent {\it Proof.} Fix a positive definite $(1,\,1)$-form $\omega$ on $V$ and consider the set \begin{eqnarray*}\widetilde{\cal K}_\omega:=\bigg\{\theta\in{\cal C}\,\mid\,\theta\wedge\omega_{n-p} = dV_\omega\bigg\}\subset(\Lambda^{p,\,p}V^\star)_\R,\end{eqnarray*} where $dV_\omega:=\omega_n$ and $\omega_l:=\omega^l/l!$ for any $l$. The set $\widetilde{\cal K}_\omega$ is closed and bounded, hence compact, in the finite-dimensional vector space $(\Lambda^{p,\,p}V^\star)_\R$. For a proof of this fact, see the proof of Proposition \ref{Prop:weak-pos_nortions-hierarchy} below (where the analogue of $\widetilde{\cal K}_\omega$ in bidegree $(n-p,\,n-p)$ is denoted by ${\cal K}_\omega\cap{\cal C}$). Moreover, we have: \begin{eqnarray*}{\cal C} = \bigg\{t\,\theta\,\mid\,t\geq 0, \hspace{1ex} \theta\in\widetilde{\cal K}_\omega\bigg\}.\end{eqnarray*}

Now, let $\Omega_\nu = t_\nu\theta_\nu\in{\cal C}$, with $\nu\geq 1$, $t_\nu\geq 0$, $\theta_\nu\in\widetilde{\cal K}_\omega$, be a sequence converging to some form $\Omega\in(\Lambda^{p,\,p}V^\star)_\R$ as $\nu\to\infty$. We have to show that $\Omega\in{\cal C}$. By continuity of the wedge product, we get the following convergence: \begin{eqnarray*}t_\nu = \frac{\Omega_\nu\wedge\omega_{n-p}}{dV_\omega}\longrightarrow\frac{\Omega\wedge\omega_{n-p}}{dV_\omega}:=t \hspace{5ex} \mbox{as}\hspace{1ex}\nu\to\infty.\end{eqnarray*} Since $t\geq 0$, there are two possible cases.

\vspace{1ex}

{\it Case $1$.} Suppose $t>0$. Then \begin{eqnarray*}\theta_\nu = \frac{1}{t_\nu}\,\Omega_\nu\longrightarrow\frac{1}{t}\,\Omega:=\theta \hspace{5ex} \mbox{as}\hspace{1ex}\nu\to\infty.\end{eqnarray*} Since $\theta_\nu\in\widetilde{\cal K}_\omega$ for every $\nu\geq 1$ and since $\widetilde{\cal K}_\omega$ is closed, we infer that $\theta\in\widetilde{\cal K}_\omega$. This implies that $\Omega = t\,\theta\in{\cal C}$.

\vspace{1ex}

{\it Case $2$.} Suppose $t=0$. Then $\Omega_\nu = t_\nu\theta_\nu\longrightarrow 0$ as $\nu\to\infty$. Indeed, $t_\nu\longrightarrow 0 (=t)$ and the sequence $(\theta_\nu)_{\nu\geq 1}$ is bounded because it is contained in the bounded set $\widetilde{\cal K}_\omega$. This proves that $\Omega = 0\in{\cal C}$. \hfill $\Box$

\vspace{2ex}

Now, fix an arbitrary $\C$-basis $A=\{\alpha_1,\dots ,\alpha_n\}$ of $V^\star$. Set \begin{eqnarray*}\Lambda_I:=i\alpha_{i_1}\wedge\bar\alpha_{i_1}\wedge\dots\wedge i\alpha_{i_p}\wedge\bar\alpha_{i_p}, \hspace{5ex} I=(1\leq i_1<\dots < i_p\leq n).\end{eqnarray*} Let $W_p(A):=\C\cdot(\Lambda_I)_{|I|=p}\subset\Lambda^{p,\,p}V^\star$ and $(W_p(A))_\R:=\R\cdot(\Lambda_I)_{|I|=p}\subset(\Lambda^{p,\,p}V^\star)_\R$ be the complex, respectively real, vector subspaces generated by the $(p,\,p)$-forms $\Lambda_I$. 

\begin{Def}\label{Def:cone_0} Consider the following set: \begin{eqnarray}\label{eqn:cone_0}{\cal C}_0(A):=\bigg\{\sum\limits_{|I|=p} a_I\Lambda_I \,\mid\,a_I\geq 0\bigg\}\subset(W_p(A))_\R\subset(\Lambda^{p,\,p}V^\star)_\R.\end{eqnarray}

\end{Def}

This is the simplicial (closed convex) cone in $(W_p(A))_\R$ generated by the $N$ rays $\R_{\geq 0}\,\Lambda_I$.

\vspace{2ex}

On the other hand, thanks to Lemma 1.4. in chapter III of Demailly's book [Dem97], there exists a $\C$-basis $\widetilde\gamma=(\gamma_s)_{1\leq s\leq N^2}$ of $\Lambda^{p,\,p}V^\star$ consisting of strongly semi-positive $(p,\,p)$-forms: \begin{eqnarray}\label{eqn:gamma_s_def}\gamma_s:=i\,\gamma_{s,\,1}\wedge\bar\gamma_{s,\,1}\wedge\dots\wedge i\,\gamma_{s,\,p}\wedge\bar\gamma_{s,\,p},\end{eqnarray} where each $\gamma_{s,\,l}$ is of one of the following shapes: $\alpha_j\pm\alpha_k$ or $\alpha_j\pm i\,\alpha_k$ for some $ 1\leq j,k\leq n$. (Taking all the possible $(1,\,0)$-forms $\alpha_j\pm\alpha_k$ and $\alpha_j\pm i\,\alpha_k$ for all $ 1\leq j,k\leq n$ as $\gamma_{s,\,l}$'s and all the $(p,\,p)$-forms $\gamma_s$ induced by them produces a system of generators of $\Lambda^{p,\,p}V^\star$, as shown in [Dem97]. So, one extracts a basis $\widetilde\gamma=(\gamma_s)_{1\leq s\leq N^2}$ from these generators. The choice of such a basis $\widetilde\gamma$ is, of course, non-unique and non-canonical. One just makes an arbitrary choice.) Since all the forms $\gamma_s$ are real and $\mbox{dim}_\R(\Lambda^{p,\,p}V^\star)_\R = N^2$, the set $\widetilde\gamma = (\gamma_s)_{1\leq s\leq N^2}$ is also an $\R$-basis of $(\Lambda^{p,\,p}V^\star)_\R$.

\begin{Def}\label{Def:cone_basis} Consider the following set: \begin{eqnarray}\label{eqn:cone_basis}{\cal C}_{A,\,\widetilde\gamma}:=\bigg\{\sum\limits_{s=1}^{N^2} a_s\gamma_s \,\mid\,a_s\geq 0\bigg\}\subset(\Lambda^{p,\,p}V^\star)_\R.\end{eqnarray}

\end{Def}  

This is a simplicial (closed convex) cone in $(\Lambda^{p,\,p}V^\star)_\R$. The definitions of ${\cal C}_0(A)$ and ${\cal C}_{A,\,\widetilde\gamma}$ suggest that they depend on the basis $A=\{\alpha_1,\dots ,\alpha_n\}$ of $V^\star$ and, in the case of the latter cone, also on the basis $\widetilde\gamma$ of $\Lambda^{p,\,p}V^\star$. We now investigate these dependencies.

A standard fact in convex geometry about simplicial cones ensures that their interiors, $\mathring{\cal C}_0(A)$ and $\mathring{\cal C}_{A,\,\widetilde\gamma}$, taken relative to their respective linear spans, $(W_p(A))_\R$ and $(\Lambda^{p,\,p}V^\star)_\R$, consist precisely of the linear combinations of all their linearly independent generators with positive coefficients: \begin{eqnarray*}\label{eqn:cone_0-basis_interior}\mathring{\cal C}_0(A) & = & \bigg\{\sum\limits_{|I|=p} a_I\Lambda_I \,\mid\, a_I>0\bigg\}\subset(W_p(A))_\R\subset(\Lambda^{p,\,p}V^\star)_\R, \\
  \mathring{\cal C}_{A,\,\widetilde\gamma} & = & \bigg\{\sum\limits_{s=1}^{N^2} a_s\gamma_s \,\mid\,a_s> 0\bigg\}\subset(\Lambda^{p,\,p}V^\star)_\R.\end{eqnarray*} Note that the interior of ${\cal C}_0(A)$ in $(\Lambda^{p,\,p}V^\star)_\R$ is empty because ${\cal C}_0(A)$ is contained in the strict vector subspace $(W_p(A))_\R\subset(\Lambda^{p,\,p}V^\star)_\R$.

\begin{Lem-Def}\label{Lem-Def:cones_identical} For every $\C$-basis $A=\{\alpha_1,\dots ,\alpha_n\}$ of $V^\star$ and every $\C$-basis $\widetilde\gamma=(\gamma_s)_{1\leq s\leq N^2}$ of $\Lambda^{p,\,p}V^\star$, the following inclusions of closed convex cones in $(\Lambda^{p,\,p}V^\star)_\R$ hold: \begin{eqnarray*}{\cal C}_0(A) \subset {\cal C} \hspace{3ex}\mbox{and}\hspace{3ex} {\cal C}_{A,\,\widetilde\gamma} \subset {\cal C}.\end{eqnarray*}

  Consequently, the interiors (in $(\Lambda^{p,\,p}V^\star)_\R$) of the last two cones satisfy the inclusion: \begin{eqnarray*}\mathring{\cal C}_{A,\,\widetilde\gamma} & = & \bigg\{\sum\limits_{s=1}^{N^2} a_s\gamma_s \,\mid\,a_s> 0\bigg\}\subset \mathring{\cal C} \subset(\Lambda^{p,\,p}V^\star)_\R.\end{eqnarray*} Moreover, the interior (in $(W_p(A))_\R$) of the first cone satisfies the inclusion: \begin{eqnarray*}\mathring{\cal C}_0(A) \subset  \mathring{\cal C}\cap(W_p(A))_\R.\end{eqnarray*}

    \vspace{1ex}

 The open convex cone $\mathring{\cal C}$ in $(\Lambda^{p,\,p}V^\star)_\R$ is independent of the choice of $\C$-basis $A=\{\alpha_1,\dots ,\alpha_n\}$ of $V^\star$. It is called the {\bf cone of strongly positive $(p,\,p)$-forms} on $V$.

\end{Lem-Def}  

\noindent {\it Proof.} The inclusion ${\cal C}_0(A)\subset{\cal C}$ follows at once from the fact that each $\Lambda_I$ is decomposable.

Similarly, the inclusion ${\cal C}_{A,\,\widetilde\gamma}\subset{\cal C}$ follows at once from the fact that each $\gamma_s$ is decomposable.  \hfill $\Box$

\vspace{2ex}

\begin{Rem}\label{Rem:reinterpretation} In the language of the intrinsic cone ${\cal C}$, Definition \ref{Def:strong-strict-positivity_pp-forms_v-space} below can be reworded as: \begin{eqnarray*}\Omega \hspace{2ex}\mbox{is {\bf strongly positive}} & \iff & \exists\,\omega\hspace{1ex}(1,\,1)\mbox{-form with}\hspace{1ex}\omega>0, \hspace{1ex} \exists\,\varepsilon>0 \hspace{2ex}\mbox{such that}\hspace{2ex} \Omega - \varepsilon\,\omega^p\in{\cal C}\\
    & \iff & \Omega\in\mathring{\cal C}.\end{eqnarray*}

\end{Rem}

\vspace{2ex}

The cone ${\cal C}_0(A)$ depends on the choice of $\C$-basis $A=\{\alpha_1,\dots ,\alpha_n\}$ of $V^\star$. However, we shall now see that the cones ${\cal C}_0(A)$ induced by all possible choices of bases $A$ exhaust the cone ${\cal C}$.

\begin{Obs}\label{Obs:union_C_0_C} The following equality holds: \begin{eqnarray}\label{eqn:union_C_0_C}\mbox{Cone}\bigg(\bigcup\limits_A{\cal C}_0(A)\bigg) = {\cal C},\end{eqnarray} where the union is taken over all $\C$-bases $A=\{\alpha_1,\dots ,\alpha_n\}$ of $V^\star$ and the left-hand side denotes the convex cone generated in $(\Lambda^{p,\,p}V^\star)_\R$ by this union of all coordinate slice cones ${\cal C}_0(A)$.

\end{Obs}

\noindent {\it Proof.} The inclusion ${\cal C}_0(A)\subset{\cal C}$ for every $\C$-basis $A$ of $V^\star$ was proved in Lemma and Definition \ref{Lem-Def:cones_identical}. So, we are reduced to proving the inclusion ``$\supset$`` in (\ref{eqn:union_C_0_C}).

Since ${\cal C}$ is generated by the decomposable forms $\Xi = i\beta_1\wedge\bar\beta_1\wedge\dots\wedge i\beta_p\wedge\bar\beta_p$, it suffices to prove that each such generator lies in ${\cal C}_0(A)$ for a suitable basis $A$. 

Let us fix a generator of ${\cal C}$ as above, with $\beta_1,\dots , \beta_p\in V^\star$. We may assume that $\Xi\neq 0$, otherwise there is nothing to prove. Then, the $(1,\,0)$-forms $\beta_1,\dots , \beta_p$ are linearly independent over $\C$, so we can complete the set they form to a $\C$-basis $B=\{\beta_1,\dots , \beta_p,\,\beta_{p+1},\dots , \beta_n\}$ of $V^\star$. Thus, we get: \begin{eqnarray*}\Xi = i\beta_1\wedge\bar\beta_1\wedge\dots\wedge i\beta_p\wedge\bar\beta_p\in{\cal C}_0(B)\end{eqnarray*} by the very definition of this cone (cf. (\ref{Def:cone_0})): \begin{eqnarray*}{\cal C}_0(B):=\bigg\{\sum\limits_{|I|=p} b_I\,i\beta_{i_1}\wedge\bar\beta_{i_1}\wedge\dots\wedge i\beta_{i_p}\wedge\bar\beta_{i_p} \,\mid\,b_I\geq 0\bigg\}.\end{eqnarray*}

This proves the inclusion ${\cal C}\subset\mbox{Cone}\bigg(\cup_A{\cal C}_0(A)\bigg)$ and completes the proof. \hfill $\Box$

\subsubsection{Our definitions for forms: (strict) positivity}\label{subsubsection:our-def}

\hspace{2ex} (I)\, First, we define {\bf strong (strict) positivity} for $C^\infty$ $(p,\,p)$-forms on a complex manifold. The same definition makes sense at a single point as well, so we start with the abstract setting of a $\C$-vector space. (The adjective/adverb ``strict/strictly'' is not necessary when writing for an English-speaking readership, so we can suppress it, but we may occasionally feel the need to include it for emphasis and for a French-speaking readership who might be influenced by the Lelong-Demailly French terminology.)

\begin{Def}\label{Def:strong-strict-positivity_pp-forms_v-space} Let $V$ be an $n$-dimensional $\C$-vector space and let $p\in\{1,\dots , n\}$. Fix a positive definite $(1,\,1)$-form $\omega\in\Lambda^{1,\,1}V^\star$ on $V$.

  We say that a $(p,\,p)$-form $\Omega\in\Lambda^{p,\,p}V^\star$ is {\bf strongly (strictly) positive} if there exists a constant $\varepsilon>0$ such that $\Omega - \varepsilon\,\omega^p$ is strongly semi-positive as a $(p,\,p)$-form on $V$.

  We write in this case: $\Omega>0$ (strongly) or $\Omega - \varepsilon\,\omega^p\geq 0$ (strongly).
  
\end{Def}

Note that, due to the fact that any two positive definite $(1,\,1)$-forms on $V$ are comparable (via positive constants), the above definition is independent of the choice of positive definite  $\omega\in\Lambda^{1,\,1}V^\star$.

\vspace{2ex}

The manifold version of this notion is introduced in the following

\begin{Def}\label{Def:strong-strict-positivity_pp-forms_manifold} Let $X$ be an $n$-dimensional complex manifold and let $\Omega\in C^\infty_{p,\,p}(X,\,\R)$ for some $p\in\{1,\dots , n\}$.

  \vspace{1ex}

  $(1)$\, Fix a point $x\in X$ and a positive definite $(1,\,1)$-form $\omega\in\Lambda^{1,\,1}T_xX^\star$ on $T^{1,\,0}_xX$.

  We say that $\Omega$ is {\bf strongly (strictly) positive at $x$} if there exists a constant $\varepsilon>0$ such that $\Omega(x) - \varepsilon\,\omega^p$ is strongly semi-positive as a $(p,\,p)$-form on $T^{1,\,0}_xX$.

  We write in this case: $\Omega(x)>0$ (strongly) or $\Omega(x) - \varepsilon\,\omega^p\geq 0$ (strongly).

\vspace{1ex}  

$(2)$\, Suppose that $X$ is compact. Let $\omega$ be a Hermitian metric on $X$. 

We say that $\Omega$ is {\bf strongly (strictly) positive} if there exists a constant $\varepsilon>0$ such that the $C^\infty$ $(p,\,p)$-form $\Omega - \varepsilon\,\omega^p$ is strongly semi-positive on $X$ (i.e. at every point of $X$).

We write in this case: $\Omega>0$ (strongly) or $\Omega - \varepsilon\,\omega^p\geq 0$ (strongly) on $X$.

\vspace{1ex}  

$(3)$\, Let $X$ be arbitrary (i.e. not necessarily compact). Let $\omega$ be a Hermitian metric on $X$.

We say that $\Omega$ is {\bf strongly (strictly) positive} if for every compact subset $K\subset X$ contained in the support of $\Omega$ there exists a constant $\varepsilon_K>0$ such that the restriction to $K$ of the $C^\infty$ $(p,\,p)$-form $\Omega - \varepsilon_K\,\omega^p$ is strongly semi-positive on $K$, namely \begin{eqnarray*}\Omega_{|K} - \varepsilon_K\,\omega^p_{|K}\geq 0 \hspace{2ex} \mbox{(strongly)} \hspace{2ex} \mbox{at every point of}\hspace{1ex} K.\end{eqnarray*}

We write in this case: $\Omega>0$ (strongly).

\end{Def}

Note that, due to the fact that any two $C^\infty$ (or merely continuous) positive definite $(1,\,1)$-forms are comparable (via positive constants) on any compact (manifold or subset of a manifold), parts $(2)$ and $(3)$ of the above definition are independent of the choice of Hermitian metric $\omega$ on $X$. Part $(1)$ too is independent of the choice of $\omega$ as it is the special case of Definition \ref{Def:strong-strict-positivity_pp-forms_v-space} for $V=T^{1,\,0}_xX$.

\begin{Rem}\label{Rem:strong-strict-positivity_pp-forms_manifold} In the setting of Definition \ref{Def:strong-strict-positivity_pp-forms_manifold}, we have the equivalences:

\vspace{1ex}  

$(1)$\, $\Omega$ is {\bf strongly (strictly) positive at $x$} $\iff$ $\Omega(x)\in\Lambda^{p,\,p}T_x^\star X$ is {\bf strongly (strictly) positive} in the sense of Definition \ref{Def:strong-strict-positivity_pp-forms_v-space} with $V=T^{1,\,0}_xX$;

\vspace{1ex}  

$(2)$\, $\Omega$ is {\bf strongly (strictly) positive} on a {\bf compact} $X$ $\iff$ for every $x\in X$, $\Omega(x)\in\Lambda^{p,\,p}T_x^\star X$ is {\bf strongly (strictly) positive} in the sense of Definition \ref{Def:strong-strict-positivity_pp-forms_v-space} with $V=T^{1,\,0}_xX$;

\vspace{1ex}  

$(3)$\, $\Omega$ is {\bf strongly (strictly) positive} on a {\bf possibly non-compact} $X$ $\iff$ for every compact subset $K\subset X$ contained in the support of $\Omega$ and for every $x\in K$,  $\Omega(x)\in\Lambda^{p,\,p}T_x^\star X$ is {\bf strongly (strictly) positive} in the sense of Definition \ref{Def:strong-strict-positivity_pp-forms_v-space} with $V=T^{1,\,0}_xX$.

\end{Rem}   

\vspace{2ex}

(II)\, We now turn our attention to {\bf weak (strict) positivity}.

\begin{Def}\label{Def:weak-strict-positivity_pp-forms_v-space} Let $V$ be an $n$-dimensional $\C$-vector space and let $\Omega\in\Lambda^{p,\,p}V^\star$ for some $p\in\{1,\dots , n\}$.

\vspace{1ex}  

$(1)$\, We say that $\Omega$ is {\bf dually weakly (strictly) positive} if \begin{eqnarray*}\Omega\wedge\Gamma >0\end{eqnarray*} for every $\Gamma\in\Lambda^{n-p,\,n-p}V^\star$ such that $\Gamma > 0$ (strongly).

  We write in this case: $\Omega>0$ (dually weakly).

\vspace{1ex}  

$(2)$\, Fix a positive definite $(1,\,1)$-form $\omega\in\Lambda^{1,\,1}V^\star$ on $V$. We say that $\Omega$ is {\bf metrically weakly (strictly) positive} if there exists a constant $\varepsilon>0$ such that $\Omega - \varepsilon\,\omega^p$ is weakly semi-positive as a $(p,\,p)$-form on $V$.

  We write in this case: $\Omega>0$ (metrically weakly) or $\Omega - \varepsilon\,\omega^p\geq 0$ (weakly).

\end{Def}

As with Definition \ref{Def:strong-strict-positivity_pp-forms_v-space}, part $(2)$ of Definition \ref{Def:weak-strict-positivity_pp-forms_v-space} is independent of the choice of positive definite $\omega\in\Lambda^{1,\,1}V^\star$. We now observe that these two notions are non-equivalent, but one implication holds.

\begin{Prop}\label{Prop:met-weakly-pos_implies_dual-weakly-pos} Let $V$ be an $n$-dimensional $\C$-vector space and let $p\in\{1,\dots , n\}$.

\vspace{1ex}
  
$(i)$\, For every $\Omega\in\Lambda^{p,\,p}V^\star$, the following implication holds: \begin{eqnarray*}\Omega \hspace{2ex} \mbox{is {\bf metrically} weakly positive} \hspace{1ex}\implies\hspace{1ex} \Omega \hspace{2ex} \mbox{is {\bf dually} weakly positive}.\end{eqnarray*}

\vspace{1ex}
  
$(ii)$\, Let $\{\alpha_1,\dots ,\alpha_n\}$ be a $\C$-basis of $V^\star$ and let $I=(1\leq i_1<\dots < i_p\leq n)$ be a fixed multi-index of length $|I|=p$. Set $\alpha_I:=\alpha_{i_1}\wedge\dots\wedge\alpha_{i_p}\in\Lambda^{p,\,0}V^\star$.

Then, the $(p,\,p)$-form \begin{eqnarray*}\Omega: = i^{p^2}\,\alpha_I\wedge\bar\alpha_I = i\alpha_{i_1}\wedge\bar\alpha_{i_1}\wedge\dots\wedge i\alpha_{i_p}\wedge\bar\alpha_{i_p}\in\Lambda^{p,\,p}V^\star\end{eqnarray*} is {\bf dually weakly positive}, but is {\bf not metrically weakly positive}.

\end{Prop}

\noindent {\it Proof.} (i)\, Suppose that $\Omega$ is metrically weakly positive. This means that $\Omega - \varepsilon\,\omega^p\geq 0$ (weakly) for some constant $\varepsilon>0$ and some positive definite $(1,\,1)$-form $\omega$. Let $\Gamma\in\Lambda^{n-p,\,n-p}V^\star$ such that $\Gamma > 0$ (strongly). Thus, there exists a constant $\delta>0$ such that $\Gamma - \delta\,\omega^{n-p}\geq 0$ (strongly).

Then, the following implication holds and the statement on the left holds as well (by the duality under (\ref{eqn:duality-pos-pairing}) between weakly semi-positive forms and strongly semi-positive forms of complementary bidegrees and the fact that this $\Gamma$ is, in particular, strongly semi-positive): \begin{eqnarray*}(\Omega - \varepsilon\,\omega^p)\wedge\Gamma\geq 0  \implies \Omega\wedge\Gamma\geq (\varepsilon\delta)\,\omega^n>0.\end{eqnarray*} To get the implication, we also used the fact that $\Gamma \geq \delta\,\omega^{n-p}$ (strongly) implies $\omega^p\wedge\Gamma\geq\delta\,\omega^n$ (strongly $=$ weakly, since these are $(n,\,n)$-forms) because any product of strongly semi-positive forms is strongly semi-positive.

Since $\Gamma$ was arbitrary with those properties, we conclude that $\Omega$ is dually weakly positive.

\vspace{1ex}

(ii)\, Fix the positive definite $(1,\,1)$-form $\omega:=\sum_{j=1}^n i\alpha_j\wedge\bar\alpha_j >0$. Then $\omega_{n-p}:=\omega^{n-p}/(n-p)! = \sum_{|J|=n-p}i^{(n-p)^2}\,\alpha_J\wedge\bar\alpha_J$.

Let $\Gamma\in\Lambda^{n-p,\,n-p}V^\star$ such that $\Gamma > 0$ (strongly). Thus, there exists a constant $\varepsilon>0$ such that $\Gamma - \varepsilon\,\omega_{n-p} \geq 0$ (strongly). Since $\Omega = i^{p^2}\,\alpha_I\wedge\bar\alpha_I\geq 0$ (strongly) -- $\Omega$ is even decomposable, we get the following inequality: \begin{eqnarray*}0\leq\Omega\wedge(\Gamma - \varepsilon\,\omega_{n-p}) = \Omega\wedge\Gamma - \varepsilon\,\Omega\wedge\omega_{n-p}.\end{eqnarray*} Hence \begin{eqnarray*}\Omega\wedge\Gamma \geq \varepsilon\,\Omega\wedge\omega_{n-p} = \varepsilon\,\bigg(i^{p^2}\,\alpha_I\wedge\bar\alpha_I\bigg)\wedge\bigg(i^{(n-p)^2}\,\alpha_{C_I}\wedge\bar\alpha_{C_I}\bigg) = \varepsilon\,dV_n >0,\end{eqnarray*} where $C_I = \{1,\dots , n\}\setminus I$ is the multi-index (of length $n-p$) complementary to $I$ and $dV_n:=i\alpha_1\wedge\bar\alpha_1\wedge\dots\wedge i\alpha_n\wedge\bar\alpha_n>0$ is the volume form induced on $V$ by the given basis of $V^\star$. This proves that $\Omega$ is {\bf dually weakly positive}.

\vspace{1ex}

If $\Omega$ were metrically weakly positive, there would exist a constant $\delta>0$ such that $\Omega - \delta\,\omega^p \geq 0$ (weakly). This would amount to \begin{eqnarray*}i^{p^2}\,\alpha_I\wedge\bar\alpha_I - \delta\, \sum_{|K|=p}i^{p^2}\,\alpha_K\wedge\bar\alpha_K \geq 0  \hspace{5ex} \mbox{(weakly)},\end{eqnarray*} which is false. Indeed, for any multi-index $L\neq I$ with $|L|=p$, we have \begin{eqnarray*}\bigg(i^{p^2}\,\alpha_I\wedge\bar\alpha_I - \delta\, \sum_{|K|=p}i^{p^2}\,\alpha_K\wedge\bar\alpha_K\bigg)\wedge\bigg(i^{(n-p)^2}\,\alpha_{C_L}\wedge\bar\alpha_{C_L}\bigg) = -\delta\,dV_n < 0.\end{eqnarray*} We conclude that $\Omega$ is {\bf not metrically weakly positive}.  \hfill $\Box$

\vspace{2ex}

We now observe a characterisation of dually weak positivity.

\begin{Prop}\label{Prop:dually-weak_char} Let $V$ be an $n$-dimensional $\C$-vector space and let $p\in\{1,\dots , n\}$. Fix a positive definite $(1,\,1)$-form $\omega\in\Lambda^{1,\,1}V^\star$ on $V$.

  For any $\Omega\in\Lambda^{p,\,p}V^\star$, the following equivalence holds: \begin{eqnarray*}\Omega > 0 \hspace{2ex} \mbox{(dually weakly)} \hspace{2ex} \iff \Omega\wedge\omega^{n-p}>0 \hspace{2ex} \mbox{and} \hspace{2ex} \Omega \geq 0 \hspace{2ex} \mbox{(weakly)}.\end{eqnarray*}

\end{Prop}  

\noindent {\it Proof.} ``$\implies$'' Suppose that $\Omega > 0$ (dually weakly). By definition, $\Omega\wedge\Gamma>0$ for every $\Gamma\in\Lambda^{n-p,\,n-p}V^\star$ such that $\Gamma>0$ (strongly). Since $\omega^{n-p}>0$ (strongly), we get $\Omega\wedge\omega^{n-p}>0$.

Showing that $\Omega\geq 0$ (weakly) is equivalent, thanks to the duality given by the pairing (\ref{eqn:duality-pos-pairing}), to showing that $\Omega\wedge v\geq 0$ for all {\it strongly semi-positive} forms $v\in\Lambda^{n-p,\,n-p}V^\star$. It then suffices to show that $\Omega$ pairs non-negatively with every decomposable $(n-p,\,n-p)$-form, namely that \begin{eqnarray}\label{eqn:dually-weak_char_proof_1}\Omega\wedge i^{(n-p)^2}\,dz_K\wedge d\bar{z}_K\geq 0\end{eqnarray} for every multi-index $K\subset\{1,\dots , n\}$ with $|K|=n-p$ and every $\C$-basis basis $\{dz_1,\dots , dz_n\}$ of $V^\star$.

  Fix such a $K$ and such a basis. Then, for every constant $1>\varepsilon>0$, we have: \begin{eqnarray*}\Omega\wedge\bigg(i^{(n-p)^2}\,dz_K\wedge d\bar{z}_K + \varepsilon\,\sum\limits_{|L|=n-p; \,L\neq K} i^{(n-p)^2}\,dz_L\wedge d\bar{z}_L\bigg)> 0\end{eqnarray*} because the $(n-p,\,n-p)$-form in the above parenthesis is strongly positive since it is ``$\geq 0$ (strongly)'' than $\varepsilon\,\gamma^{n-p}/(n-p)!$, where $\gamma:=\sum_{1\leq j\leq n} idz_j\wedge d\bar{z}_j>0$. Letting $\varepsilon\searrow 0$, we get (\ref{eqn:dually-weak_char_proof_1}).

\vspace{1ex}

``$\Longleftarrow$'' Let $\Gamma\in\Lambda^{n-p,\,n-p}V^\star$ such that $\Gamma>0$ (strongly). Then, there exists a constant $\varepsilon>0$ and a form $\Gamma_\varepsilon\in\Lambda^{n-p,\,n-p}V^\star$ such that $\Gamma_\varepsilon\geq 0$ (strongly) and $\Gamma = \varepsilon\,\omega^{n-p} + \Gamma_\varepsilon$. We get: \begin{eqnarray*}\Omega\wedge\Gamma = \varepsilon\,\Omega\wedge\omega^{n-p} + \Omega\wedge\Gamma_\varepsilon >0\end{eqnarray*} because $\Omega\wedge\omega^{n-p}>0$ (by hypothesis) and $\Omega\wedge\Gamma_\varepsilon\geq 0$ (by hypothesis and the strong semi-positivity of $\Gamma_\varepsilon$). We conclude that $\Omega\wedge\Gamma>0$ for every $\Gamma\in\Lambda^{n-p,\,n-p}V^\star$ such that $\Gamma>0$ (strongly), which amounts to $\Omega > 0$ (strongly). \hfill $\Box$

\vspace{2ex}

We now introduce a third notion of ``weak (strict) positivity'' defined by duality with {\bf decomposable} forms of the complementary bidegree. (Cf. the {\it dually weak positivity} of Definition \ref{Def:weak-strict-positivity_pp-forms_v-space} that is defined by duality with {\it strongly positive} forms  of the complementary bidegree.) 

\begin{Def}\label{Def:decomposably-weak-strict-positivity_pp-forms_v-space} Let $V$ be an $n$-dimensional $\C$-vector space and let $\Omega\in\Lambda^{p,\,p}V^\star$ for some $p\in\{1,\dots , n\}$. We say that $\Omega$ is {\bf decomposably weakly (strictly) positive} if for any $\C$-linearly independent $(1,\,0)$-forms $\alpha_1,\dots , \alpha_{n-p}\in V^\star$, we have \begin{eqnarray*}\Omega\wedge i\alpha_1\wedge\bar\alpha_1\wedge\dots\wedge i\alpha_{n-p}\wedge\bar\alpha_{n-p} >0.\end{eqnarray*}

  We write in this case: $\Omega>0$ (decomposably weakly).

\end{Def}

Note that it suffices to check the property $\Omega\wedge i^{(n-p)^2}\,\beta_K\wedge\bar\beta_K>0$ for every multi-index $K$ with $|K|=n-p$ and for a {\it fixed} $\C$-basis $\{\beta_1,\dots , \beta_n\}$ of $V^\star$.

This notion has the following geometric interpretation, analogous to the classical one for {\it weakly semi-positive} forms.

\begin{Prop}\label{Prop:decomposably-weak-strict-positivity_restrictions} In the setting of Definition \ref{Def:decomposably-weak-strict-positivity_pp-forms_v-space}, for any $\Omega\in\Lambda^{p,\,p}V^\star$ the next equivalence holds: \begin{eqnarray*}\Omega > 0 \hspace{2ex} \mbox{(decomposably weakly)} & \iff & \Omega_{|S}>0 \hspace{2ex} \mbox{as a volume form on} \hspace{2ex} S, \\
    & & \hspace{9ex} \mbox{for every $p$-dimensional $\C$-vector subspace} \hspace{2ex} S\subset V.\end{eqnarray*}  

\end{Prop}

\noindent {\it Proof.}  ``$\implies$'' Suppose that $\Omega > 0$ (decomposably weakly). Let $S\subset V$ be a $p$-dimensional $\C$-vector subspace. (So, $\mbox{codim}_V S = n-p$.) There exists a $\C$-basis $\{\alpha_1,\dots , \alpha_{n-p},\, \alpha_{n-p+1},\dots , \alpha_n\}$ of $V^\star$ such that $S = \{\alpha_1 = \dots = \alpha_{n-p} =0\}\subset V$. We write \begin{eqnarray*}\Omega = \sum\limits_{|I| = |J| = p}\Omega_{I\bar{J}}\,i^{p^2}\,\alpha_I\wedge\bar\alpha_J,   \hspace{5ex} \Omega_{I\bar{J}}\in\C.\end{eqnarray*}

Then, on the one hand, the following general identity holds: \begin{eqnarray}\label{eqn:restriction_S_coeff-form}\Omega_{|S} = \Omega_{L\bar{L}}\,i^{p^2}\,\alpha_L\wedge\bar\alpha_L, \hspace{5ex} \mbox{where}\hspace{2ex} L= (n-p+1<\dots < n).\end{eqnarray} On the other hand, the following general identity holds: \begin{eqnarray}\label{eqn:Omega_wedge_decomposable_S}\Omega\wedge i\alpha_1\wedge\bar\alpha_1\wedge\dots\wedge i\alpha_{n-p}\wedge\bar\alpha_{n-p} = \Omega_{L\bar{L}}\,dV_n, \hspace{5ex} \mbox{where}\hspace{2ex} L= (n-p+1<\dots < n).\end{eqnarray} We have set $dV_n=i\alpha_1\wedge\bar\alpha_1\wedge\dots\wedge i\alpha_n\wedge\bar\alpha_n>0$, the volume form induced by the chosen basis of $V^\star$.

Now, the decomposably weak positivity hypothesis on $\Omega$ implies the positivity of the volume form on the left-hand side of (\ref{eqn:Omega_wedge_decomposable_S}). This yields $\Omega_{L\bar{L}}>0$. This, in turn, implies the positivity (as a volume form on $S$) of the right-hand side of (\ref{eqn:restriction_S_coeff-form}). Consequently, $\Omega_{|S}>0$, as desired.

\vspace{1ex}

``$\Longleftarrow$'' Let $\alpha_1,\dots , \alpha_{n-p}\in V^\star$ be $\C$-linearly independent $(1,\,0)$-forms on $V$. Complete this set to a $\C$-basis $\{\alpha_1,\dots , \alpha_{n-p},\,\alpha_{n-p+1},\dots , \alpha_n\}$ of $V^\star$. Consider the $p$-dimensional $\C$-vector subspace $S = \{\alpha_1 = \dots = \alpha_{n-p} =0\}\subset V$. The hypothesis yields $\Omega_{|S}>0$ as a volume form on $S$. Thanks to the general identity (\ref{eqn:restriction_S_coeff-form}), this is equivalent to $\Omega_{L\bar{L}}>0$, which in turn, thanks to the general identity (\ref{eqn:Omega_wedge_decomposable_S}), is equivalent to \begin{eqnarray*}\Omega\wedge i\alpha_1\wedge\bar\alpha_1\wedge\dots\wedge i\alpha_{n-p}\wedge\bar\alpha_{n-p} > 0.\end{eqnarray*}

We conclude that $\Omega > 0$ (decomposably weakly), as desired.  \hfill $\Box$

\vspace{2ex}

We now characterise of two of these {\it weak positivity} notions in terms of evaluations of the cone ${\cal C}$.

\begin{Prop}\label{Prop:weak-pos_evaluations_cone} Let $V$ be an $n$-dimensional $\C$-vector space and let $p\in\{1,\dots , n\}$. Consider the cone ${\cal C}\subset(\Lambda^{n-p,\,n-p}V^\star)_\R$ of strongly semi-positive $(n-p,\,n-p)$-forms on $V$ introduced in Definition \ref{Def:cone_all-decomposable} (with $p$ in place of $n-p$).

  Fix a form $\Omega\in\Lambda^{p,\,p}V^\star$. The following equivalences hold: \begin{eqnarray}\label{eqn:decomp-weak-pos_evaluations_cone}\Omega>0 \hspace{1ex} \mbox{({\bf decomposably} weakly)} & \iff & \Omega\wedge\Xi>0 \hspace{2ex}\mbox{for every}\hspace{1ex} \Xi\in{\cal C}\setminus\{0\} \\
  \label{eqn:dual-weak-pos_evaluations_cone}\Omega>0 \hspace{1ex} \mbox{({\bf dually} weakly)} & \iff & \Omega\wedge\Xi>0 \hspace{2ex}\mbox{for every}\hspace{1ex} \Xi\in\mathring{\cal C}.\end{eqnarray}

\end{Prop}

\noindent {\it Proof.} Since the condition ``$\Omega>0$ (decomposably weakly)'' is equivalent to $\Omega\wedge u>0$ for every non-zero decomposable $(n-p,\,n-p)$-form $u$, it is also equivalent to $\Omega\wedge\sum_j c_j\, u_j>0$ for all finite linear combinations with non-negative coefficients $c_j\geq 0$, at least one of which is positive, of non-zero decomposable $(n-p,\,n-p)$-forms $u_j$. But these forms $\Xi:= \sum_j c_j\, u_j$ are precisely (all) the elements of ${\cal C}\setminus\{0\}$, hence the first equivalence in the statement.

The second equivalence in the statement follows at once by putting together $(1)$ of Definition \ref{Def:weak-strict-positivity_pp-forms_v-space} and Remark \ref{Rem:reinterpretation}.  \hfill $\Box$

\vspace{2ex}

We now observe that {\it decomposably weak positivity} is equivalent to {\it metrically weak positivity}.

\begin{Prop}\label{Prop:weak-pos_nortions-hierarchy} Let $V$ be an $n$-dimensional $\C$-vector space and let $p\in\{1,\dots , n\}$. For every $\Omega\in\Lambda^{p,\,p}V^\star$, the following equivalence and implication hold: \begin{eqnarray*}\Omega >0 \hspace{1ex} \mbox{({\bf metrically} weakly)} \iff \Omega>0 \hspace{1ex} \mbox{({\bf decomposably} weakly)}  \implies \Omega>0 \hspace{1ex} \mbox{({\bf dually} weakly)}.\end{eqnarray*}

\end{Prop}
  
\noindent {\it Proof.} In view of (i) of Proposition \ref{Prop:met-weakly-pos_implies_dual-weakly-pos}, it suffices to prove the stated equivalence.

$\bullet$ Suppose that $\Omega>0$ (metrically weakly). Let $\alpha_1,\dots , \alpha_{n-p}\in V^\star$ be $\C$-linearly independent $(1,\,0)$-forms. Complete this set to a $\C$-basis $\{\alpha_1,\dots , \alpha_{n-p},\,\alpha_{n-p+1},\dots , \alpha_n\}$ of $V^\star$ and consider the positive definite $(1,\,1)$-form $\omega:=\sum_{1\leq j\leq n} i\alpha_j\wedge\bar\alpha_j>0$. Then, $\omega_p:=\omega^p/p! = \sum_{|I|=p} i^{p^2}\,\alpha_I\wedge\bar\alpha_I$.

The metrically weak positivity assumption on $\Omega$ yields a constant $\varepsilon>0$ such that $\Omega - \varepsilon\,\omega_p \geq 0$ (weakly). Since $i\alpha_1\wedge\bar\alpha_1\wedge\dots\wedge i\alpha_{n-p}\wedge\bar\alpha_{n-p} \geq 0$ (strongly), by duality we get: \begin{eqnarray*}(\Omega - \varepsilon\,\omega_p)\wedge i\alpha_1\wedge\bar\alpha_1\wedge\dots\wedge i\alpha_{n-p}\wedge\bar\alpha_{n-p}\geq 0,\end{eqnarray*} which amounts to the first inequality below: \begin{eqnarray*}\Omega\wedge i\alpha_1\wedge\bar\alpha_1\wedge\dots\wedge i\alpha_{n-p}\wedge\bar\alpha_{n-p}\geq \varepsilon\,\omega_p \wedge i\alpha_1\wedge\bar\alpha_1\wedge\dots\wedge i\alpha_{n-p}\wedge\bar\alpha_{n-p} = \varepsilon\,dV_n>0,\end{eqnarray*} where $dV_n=i\alpha_1\wedge\bar\alpha_1\wedge\dots\wedge i\alpha_n\wedge\bar\alpha_n>0$ is the volume form induced by the chosen basis of $V^\star$.

This proves that $\Omega>0$ (decomposably weakly).

\vspace{1ex}

$\bullet$ Suppose that $\Omega>0$ (decomposably weakly). By (\ref{eqn:decomp-weak-pos_evaluations_cone}), this means that the linear map \begin{eqnarray*}l_{\omega,\,\Omega} : (\Lambda^{n-p,\,n-p}V^\star)_\R\longrightarrow\R, \hspace{5ex} u\longmapsto\frac{\Omega\wedge u}{dV_\omega},\end{eqnarray*} is positive on ${\cal C}\setminus\{0\}$, where $\omega>0$ is any positive definite $(1,\,1)$-form on $V$ and $dV_\omega = \omega^n/n!$ is its associated volume form.

Meanwhile, $(\Lambda^{n-p,\,n-p}V^\star)_\R$ is finite-dimensional, so $l_{\omega,\,\Omega}$ is continuous, so the affine hyperplane \begin{eqnarray*}{\cal K}_\omega:=\bigg\{\theta\in(\Lambda^{n-p,\,n-p}V^\star)_\R\,\mid\,\theta\wedge\omega_p = dV_\omega\bigg\}\subset(\Lambda^{n-p,\,n-p}V^\star)_\R\end{eqnarray*} is closed. Thus, ${\cal K}_\omega\cap{\cal C}$ is closed and bounded in the finite-dimensional vector space $(\Lambda^{n-p,\,n-p}V^\star)_\R$, hence compact.

  Proving that $\Omega>0$ (metrically weakly) is equivalent to proving the existence of a constant $\varepsilon>0$ such that $\Omega\wedge u\geq \varepsilon\,\omega_p\wedge u$ for every decomposable $(n-p,\,n-p)$-form $u$ on $V$. In other words, we have the equivalence: \begin{eqnarray*}\Omega >0 \hspace{1ex} \mbox{(metrically weakly)} \iff \exists\,\varepsilon>0 \hspace{1ex}\mbox{such that}\hspace{1ex} l_{\omega,\,\Omega}\geq\varepsilon \hspace{2ex}\mbox{on all decomposables in}\hspace{1ex} {\cal K}_\omega\cap{\cal C}.\end{eqnarray*}

  Now, under our assumption, $l_{\omega,\,\Omega}>0$ on ${\cal C}\setminus\{0\}$, hence also on the compact subset ${\cal K}_\omega\cap{\cal C}$. By continuity of $l_{\omega,\,\Omega}$ and compactness of ${\cal K}_\omega\cap{\cal C}$, there exists $\varepsilon>0$ such that $l_{\omega,\,\Omega}\geq\varepsilon$ on ${\cal K}_\omega\cap{\cal C}$.

  This proves that  $\Omega>0$ (metrically weakly). The key and standard fact that ${\cal K}_\omega\cap{\cal C}$ is bounded, which has been used, can be proved by using, for example, the cone equality (\ref{eqn:union_C_0_C}). Indeed, for every $\C$-basis $A=\{\alpha_1,\dots ,\alpha_n\}$ of $V^\star$, normalised such that $\omega = \sum_{1\leq j\leq n} i\alpha_j\wedge\bar\alpha_j$, and every $\theta=\sum\limits_{|I|=n-p} a_I\Lambda_I\in{\cal C}_0(A)$, we have $a_I\geq 0$ for all $I$ and \begin{eqnarray*}\theta\wedge\omega_p = \bigg(\sum\limits_{|I|=n-p} a_I\bigg)\,dV_\omega.\end{eqnarray*} Thus, for any $\theta\in{\cal C}_0(A)$, its membership in ${\cal K}_\omega$ is equivalent to $\sum_{|I|=n-p} a_I =1$ and implies that $0\leq a_I\leq 1$ for every $I$. This proves that the ``slices'' ${\cal K}_\omega\cap{\cal C}_0(A)$ of ${\cal K}_\omega\cap{\cal C}$ are uniformly bounded, so ${\cal K}_\omega\cap{\cal C}$ is bounded.  \hfill $\Box$

\vspace{2ex}

We now observe that the relations among the various notions of strict positivity (all of which are independent of the choice of basis of $V^\star$) simplify drastically in the case of {\bf diagonalisable} forms (cf. Observation \ref{Obs:m-m_diagonal_equiv}).

\begin{Obs}\label{Obs:equiv_diagonal-forms_strong-weak} Let $V$ be an $n$-dimensional $\C$-vector space and let $p\in\{1,\dots , n\}$. Fix a $\C$-basis $\{\alpha_1,\dots , \alpha_n\}$ of $V^\star$ and a  $(p,\,p)$-form that is {\bf diagonal} with respect to this basis: \begin{eqnarray*}\Omega = \sum\limits_{|I|=p}\lambda_I\,i^{p^2}\,\alpha_I\wedge\bar\alpha_I = \sum\limits_{|I|=p}\lambda_I\,i\alpha_{i_1}\wedge\bar\alpha_{i_1}\wedge\dots\wedge i\alpha_{i_p}\wedge\bar\alpha_{i_p}\in\Lambda^{p,\,p}V^\star, \hspace{5ex} \lambda_I\in\C.\end{eqnarray*}

  (a)\, Then, the following equivalences hold: \begin{eqnarray*}\Omega> 0 \hspace{2ex} \mbox{(strongly)} & \iff & \Omega >0 \hspace{1ex} \mbox{({\bf metrically} weakly)} \iff \Omega>0 \hspace{1ex} \mbox{({\bf decomposably} weakly)} \\
    & \iff &  \lambda_I> 0 \hspace{5ex} \mbox{for every}\hspace{1ex} I  \hspace{1ex} \mbox{with}\hspace{1ex} |I|=p.\end{eqnarray*}

  (b)\, Meanwhile, the following equivalence holds: \begin{eqnarray*} \Omega> 0 \hspace{2ex} \mbox{({\bf dually} weakly)} \iff \lambda_I\geq 0 \hspace{5ex} \mbox{for every}\hspace{1ex} I  \hspace{1ex} \mbox{with}\hspace{1ex} |I|=p \hspace{2ex} \mbox{and}\hspace{2ex} \sum\limits_{|I|=p}\lambda_I>0.\end{eqnarray*}

\end{Obs}

\noindent {\it Proof.} (a)\, We have already seen the implications: \begin{eqnarray*}\Omega> 0 \hspace{2ex} \mbox{(strongly)} \implies \Omega >0 \hspace{1ex} \mbox{(metrically weakly)} \implies \Omega>0 \hspace{1ex} \mbox{(decomposably weakly)}.\end{eqnarray*}

Now, if $\Omega>0$ (decomposably weakly), for every multi-index $K$ with $|K|=n-p$ we get: \begin{eqnarray*}0<\Omega\wedge i^{(n-p)^2}\,dz_K\wedge d\bar{z}_K = \lambda_{C_K}\,dV_n,\end{eqnarray*} where $C_K$ is the multi-index (of length $p$) that is complementary to $K$. We infer that $\lambda_{C_K}>0$. Thus, $\lambda_I>0$ for every $I$ with $|I|=p$.

Finally, if $\lambda_I>0$ for every $I$ with $|I|=p$, let $\varepsilon:=\min\,\{\lambda_I\,\mid\,|I|=p\}>0$. Denoting by $\omega:=\sum_{1\leq j\leq n} i\alpha_j\wedge\bar\alpha_j>0$, we get: \begin{eqnarray*}\Omega \geq \varepsilon\,\sum\limits_{|I|=p}i^{p^2}\,\alpha_I\wedge\bar\alpha_I = \varepsilon\,\omega_p  \hspace{5ex} \mbox{(strongly)},\end{eqnarray*} proving that $\Omega>0$ (strongly). The ``$\geq$'' inequality above holds in the strong sense because each $(p,\,p)$-form $i^{p^2}\,\alpha_I\wedge\bar\alpha_I$ is decomposable, hence also strongly semi-positive.

\hspace{1ex}

(b)\, For every $I$ with $|I|=p$, the inequality $\lambda_I\geq 0$ is equivalent to $\Omega\wedge i^{(n-p)^2}\,dz_K\wedge d\bar{z}_K\geq 0$, where $K=C_I$. Thus, $\lambda_I\geq 0$ for all $I$ with $|I|=p$ if and only if $\Omega$ pairs non-negatively with all decomposables of bidegree $(n-p,\,n-p)$ induced by the basis $\{\alpha_1,\dots , \alpha_n\}$ of $V^\star$. This is equivalent to $\Omega\wedge\Gamma\geq 0$ for every decomposable $(n-p,\,n-p)$-form, hence to $\Omega\geq 0$ (weakly).

On the other hand, $\omega_{n-p} = \sum_{|J|=n-p}i^{(n-p)^2}\,\alpha_J\wedge\bar\alpha_J$, so \begin{eqnarray*}\Omega\wedge\omega_{n-p} = \bigg(\sum\limits_{|I|=p}\lambda_I\bigg)\,dV_n.\end{eqnarray*}

Putting these two pieces of information together, we get the equivalence: \begin{eqnarray*}\lambda_I\geq 0 \hspace{2ex} \mbox{for every}\hspace{1ex} I  \hspace{1ex} \mbox{with}\hspace{1ex} |I|=p \hspace{3ex} \mbox{and}\hspace{3ex} \sum\limits_{|I|=p}\lambda_I>0 \iff \Omega\geq 0 \hspace{2ex} \mbox{(weakly)} \hspace{2ex} \mbox{and} \hspace{2ex} \Omega\wedge\omega_{n-p}>0.\end{eqnarray*}

The contention follows from this equivalence and from the characterisation of {\it dually weak positivity} given in Proposition \ref{Prop:dually-weak_char}. \hfill $\Box$

\vspace{2ex}

Since in bidegrees $(1,\,1)$ and $(n-1,\,n-1)$ all forms are diagonalisable, Observation \ref{Obs:equiv_diagonal-forms_strong-weak} applies to all those forms and we get the

\begin{Cor}\label{Obs:equiv_diagonal-forms_strong-weak_11} For every form $\Omega\in\Lambda^{1,\,1}V^\star\cup\Lambda^{n-1,\,n-1}V^\star$, the following equivalences hold: \begin{eqnarray*}\Omega> 0 \hspace{2ex} \mbox{(strongly)} & \iff & \Omega >0 \hspace{1ex} \mbox{({\bf metrically} weakly)} \iff \Omega>0 \hspace{1ex} \mbox{({\bf decomposably} weakly)} \\
    & \iff &  \mbox{the coefficient matrix of $\Omega$ in some, hence any, basis of $V^\star$} \\
    & & \mbox{is {\bf positive definite}} \end{eqnarray*} and \begin{eqnarray*} \Omega> 0 \hspace{2ex} \mbox{({\bf dually} weakly)} & \iff & \mbox{all the eigenvalues of $\Omega$ w.r.t. some, hence any, metric $\omega>0$}  \\
    & & \mbox{are {\bf non-negative} and the trace of $\Omega$ is {\bf positive}}.\end{eqnarray*}

\end{Cor}

\vspace{2ex}

We also observe the following

\begin{Prop}\label{Prop:Omega-prod-omega_m-1} Let $V$ be an $n$-dimensional $\C$-vector space and let $m\in\{1,\dots , n\}$. For every $\Omega\in\Lambda^{n-m,\,n-m}V^\star$ such that $\Omega>0$ ({\bf decomposably weakly}) and every $\omega\in\Lambda^{1,\,1}V^\star$ such that $\omega>0$, we have $\Omega\wedge\omega^{m-1} > 0$ ({\bf decomposably weakly}).

    In other words, the $(n-1,\,n-1)$-form $\Omega\wedge\omega^{m-1}$ is {\bf positive definite} in the usual sense under these hypotheses.

\end{Prop}  
  
\noindent {\it Proof.} We saw in Proposition \ref{Prop:weak-pos_nortions-hierarchy} that the hypothesis ``$\Omega>0$ ( decomposably weakly)'' is equivalent to 
``$\Omega>0$ (metrically weakly)''. Thus, there exists a constant $\varepsilon>0$ such that $\Omega - \varepsilon\,\omega^{n-m}\geq 0$ (weakly). Since $\omega^{m-1} > 0$ (strongly), we get $\Omega\wedge\omega^{m-1} - \varepsilon\,\omega^{n-1}\geq 0$ (weakly), which proves that $\Omega\wedge\omega^{m-1} > 0$ (metrically weakly). This is equivalent to the contention.  \hfill $\Box$

\subsection{(Semi-)positivity for currents}\label{subsection:semi_pos_currents} We now deal with the counterparts for currents of the notions analysed in the previous subsection for forms. Let $X$ be a complex manifold with $\mbox{dim}_\C X=n$ and let $p\in\{1,\dots, n\}$. We will use the (by now) standard notation of [Dem97]:

\vspace{1ex}

-the space of compactly supported $C^\infty$ forms of bidegree $(p,\,p)$ on $X$ (also called {\it test $(p,\,p)$-forms}) is denoted by ${\cal D}_{p,\,p}(X)$; the space of compactly supported $C^\infty$ functions on $X$ (also called {\it test functions}) is denoted by ${\cal D}(X)$;

\vspace{1ex}

-for every compact subset $K\subset X$, the space of $C^\infty$ forms of bidegree $(p,\,p)$ on $X$ with support contained in $K$ is denoted by ${\cal D}_{p,\,p}(K)$; we have ${\cal D}_{p,\,p}(K)\subset{\cal D}_{p,\,p}(X)$ and ${\cal D}_{p,\,p}(X) = \cup_K{\cal D}_{p,\,p}(K)$;

\vspace{1ex}

-the space of currents of bidegree $(n-p,\,n-p)$ (or, equivalently, of bidimension $(p,\,p)$) on $X$ is denoted by ${\cal D}^{'n-p,\,n-p}(X)$ or by ${\cal D}'_{p,\,p}(X)$; the space of distributions on $X$ is denoted by ${\cal D}'(X)$.

\vspace{1ex}

Recall that currents of bidegree $(n-p,\,n-p)$ are, by definition, linear maps $T:{\cal D}_{p,\,p}(X)\longrightarrow\C$ whose restriction to ${\cal D}_{p,\,p}(K)$ for each compact subset $K\subset X$ is continuous. With respect to local holomorphic coordinates $z_1,\dots , z_n$ on $X$, every current $T\in{\cal D}^{'n-p,\,n-p}(X)$ has the shape \begin{eqnarray*}T = \sum\limits_{|I| = |J| = n-p} T_{I\bar{J}}\,i^{(n-p)^2}\,dz_I\wedge d\bar{z}_J,\end{eqnarray*} where the coefficients $T_{I\bar{J}}$ are distributions on the coordinate patch $U\subset X$. They are defined by \begin{eqnarray*}T_{I\bar{J}}(f):= \langle T,\, f\,i^{p^2}\,dz_{C_I}\wedge d\bar{z}_{C_J}\rangle = \int\limits_X T\wedge f\,i^{p^2}\,dz_{C_I}\wedge d\bar{z}_{C_J}, \hspace{5ex} f\in{\cal D}(U).\end{eqnarray*} We have used the general notation and equality: \begin{eqnarray*}\langle T,\,u\rangle = \int\limits_X T\wedge u,  \hspace{5ex} T\in {\cal D}^{'n-p,\,n-p}(X), \hspace{2ex} u\in{\cal D}_{p,\,p}(X),\end{eqnarray*} where $\langle T,\,u\rangle$ denotes the pairing of a current $T$ with a test form $u$ of the complementary bidegree (given by the definition of the current) and $T\wedge u$ denotes the current of bidegree $(n,\,n)$ obtained by exterior multiplication. Indeed, a current can always be multiplied with a smooth form since distributions can be multiplied with smooth functions, although two currents cannot, in general, be multiplied. In particular, every smooth form $u\in C^\infty_{n-p,\,n-p}(X,\,\C)$ defines a current $T\in{\cal D}^{'n-p,\,n-p}(X)$ of the same bidegree.

\subsubsection{Classical notions for currents: semi-positivity}\label{subsubsection:classical_currents} We briefly review here the main notions presented in $1.B.$ of chapter III of Demailly's book [Dem97].

According to [Dem97, III, Definition 1.13.], a current $T\in{\cal D}^{'n-p,\,n-p}(X)$ is said to be:

\vspace{1ex}

(i)\, {\bf weakly semi-positive} if $\langle T,\,u\rangle\geq 0$ for every $u\in{\cal D}_{p,\,p}(X)$ such that $u(x)\geq 0$ ({\bf strongly}) at every point $x\in X$;

(ii)\, {\bf strongly semi-positive} if $\langle T,\,u\rangle\geq 0$ for every $u\in{\cal D}_{p,\,p}(X)$ such that $u(x)\geq 0$ ({\bf weakly}) at every point $x\in X$.

\vspace{1ex}

We write ``$T\geq 0$ (weakly)'' in case (i) and ``$T\geq 0$ (strongly)'' in case (ii).

\vspace{1ex}

A fundamental example is the following: the current of integration $[Z]\in{\cal D}^{'n-p,\,n-p}(X)$ defined by every $p$-dimensional complex submanifold $Z\subset X$ as \begin{eqnarray*}\bigg\langle[Z],\,u\bigg\rangle:=\int\limits_Z u_{|Z},  \hspace{5ex} u\in{\cal D}_{p,\,p}(X),\end{eqnarray*} is {\it strongly semi-positive} since, whenever $u\geq 0$ (weakly), its restriction $u_{|Z}$ to every $Z$ as above is a non-negative volume form on $Z$.

We see that these standard semi-positivity notions for currents are defined by means of the corresponding notions for forms (recalled in $\S$\ref{subsubsection:classical}) and the {\it current/test form} duality: \begin{eqnarray}\label{eqn:duality-currents-forms}{\cal D}^{'n-p,\,n-p}(X)\times{\cal D}_{p,\,p}(X)\longrightarrow\C, \hspace{5ex} (T,\,u)\longmapsto\langle T,\,u\rangle = \int\limits_X T\wedge u.\end{eqnarray}

In particular, whenever we have a sequence of currents $T_\nu\longrightarrow T$ converging weakly as $\nu\to\infty$, the limiting current $T$ is {\it weakly semi-positive} (respectively {\it strongly semi-positive}) whenever the currents $T_\nu$ are {\it weakly semi-positive} (respectively {\it strongly semi-positive}). Thus, in every bidegree $(n-p,\,n-p)$, the set of {\it weakly semi-positive} currents and the set of {\it strongly semi-positive} currents are closed convex cones in ${\cal D}^{'n-p,\,n-p}(X)$.

We refer the reader to [Dem97, III, $1.B.$] for the fundamental properties of {\it weakly} and {\it strongly semi-positive} currents. We only recall here the following key result, a consequence of the classical Banach-Alaoglu theorem.

\begin{Prop}([Dem97, III, Proposition 1.23.])\label{Prop:semi-pos-currents_b-mass_compactness} For every Hermitian metric $\omega$ on $X$, every $p\in\{1,\dots , n\}$ and every constant $C>0$, the set \begin{eqnarray*}\bigg\{T\in{\cal D}^{'n-p,\,n-p}(X)\,\mid\, T\geq 0 \hspace{1ex} \mbox{(weakly)} \hspace{2ex}\mbox{and}\hspace{2ex} \int\limits_X T\wedge\omega^p\leq C\bigg\}\end{eqnarray*} is {\bf compact} in ${\cal D}^{'n-p,\,n-p}(X)$ with respect to the weak topology of currents.

\end{Prop}

\vspace{1ex}

Still, before moving on to introducing our notions of {\it weak} and {\it strong positivity} for currents in the next subsection, we point out that a $C^\infty$ form is {\it weakly} or {\it strongly semi-positive} as a form at every point on $X$ (in the sense recalled in $\S$\ref{subsubsection:classical} that partially uses the pointwise duality (\ref{eqn:duality-pos-pairing})) if and only if it has the same property in the sense of the integral duality (\ref{eqn:duality-currents-forms}).

\begin{Obs}\label{Obs:semi-pos_same_forms-currents} Let $\Omega\in C^\infty_{n-p,\,n-p}(X,\,\C)$. The following equivalences hold:

\vspace{1ex}

(a)\, $\Omega\geq 0$ (weakly) as a form at every $x\in X$ $\iff$ $\Omega\geq 0$ (weakly) as a current on $X$;

\vspace{1ex}

(b)\, $\Omega\geq 0$ (strongly) as a form at every $x\in X$ $\iff$ $\Omega\geq 0$ (strongly) as a current on $X$.

\end{Obs}  

\noindent {\it Proof.} We will prove (b), the proof of (a) being very similar.

\vspace{1ex}

``$\implies$'' The hypothesis implies that $\Omega(x)\wedge\Gamma(x)\geq 0$ at every point $x\in X$ for every form $\Gamma\in{\cal D}_{p,\,p}(X)$ such that $\Gamma(x) \geq 0$ (weakly) for every $x\in X$. Integrating the compactly supported smooth $(n,\,n)$-form $\Omega\wedge\Gamma$, we get:  \begin{eqnarray*}\langle\Omega,\,\Gamma\rangle = \int\limits_X\Omega\wedge\Gamma\geq 0.\end{eqnarray*} Since $\Gamma$ was arbitrary with the above properties, we infer that $\Omega\geq 0$ (strongly) as a current on $X$. 

\vspace{1ex}

``$\Longleftarrow$'' Suppose that $\Omega\geq 0$ (strongly) as a current on $X$. Fix a point $x\in X$. Reasoning by contradiction, we suppose there exists a form $\Gamma_x\in(\Lambda^{p,\,p}T_x^\star X)_\R$ such that $\Gamma_x\geq 0$ (weakly) and $\Omega(x)\wedge\Gamma_x <0$.

Let $z_1,\dots , z_n$ be local holomorphic coordinates on $X$, centred at $x$ and defined on an open subset $U\subset X$. Then \begin{eqnarray*}\Gamma_x = \sum\limits_{|J|=|K|=p}\Gamma_{J\overline{K}}\,i^{p^2}\,dz_J\wedge d\overline{z}_K,    \hspace{5ex} \Gamma_{J\bar{K}}\in\C.\end{eqnarray*} We extend $\Gamma_x$ to a form $\Gamma\in C^\infty_{p,\,p}(U,\,\R)$ by keeping its coefficients $\Gamma_{J\overline{K}}$ constant on $U$. Then $\Gamma(z)\geq 0$ (weakly) for every $z\in U$ and, by continuity, $\Omega(z)\wedge\Gamma(z) <0$ for every $z\in U_0$ for some open subset $U_0\subset U$. Fix a non-empty open subset $U_1\Subset U_0$ and a $C^\infty$ function $\psi:X\longrightarrow[0,\,1]$ with compact support such that $\psi\equiv 1$ on $U_1$ and $\mbox{Supp}\,\psi\Subset U_0$. Thus, $\psi\Gamma\in{\cal D}_{p,\,p}(X)$ satisfies $\psi\Gamma\geq 0$ (weakly) everywhere on $X$ and $\Omega\wedge(\psi\Gamma)\in{\cal D}_{n,\,n}(X)$ satisfies: \begin{eqnarray*}\Omega\wedge(\psi\Gamma) = 0 \hspace{1ex} \mbox{on}\hspace{1ex} X\setminus U_0; \hspace{3ex} \Omega\wedge(\psi\Gamma) < 0 \hspace{1ex} \mbox{on}\hspace{1ex} U_1; \hspace{3ex} \Omega\wedge(\psi\Gamma) \leq 0 \hspace{1ex} \mbox{on}\hspace{1ex} U_0.\end{eqnarray*}

We infer that $\int_X\Omega\wedge(\psi\Gamma) < 0$, contradicting the hypothesis $\Omega\geq 0$ (strongly) as a current on $X$.  \hfill $\Box$

\subsubsection{Our definitions for currents: (strict) positivity}\label{subsubsection:our-def_currents}

Following the definitions given in $\S$\ref{subsubsection:our-def} for forms, we now define the analogous notions for currents.

\begin{Def}\label{Def:strict-positivity_currents} Let $X$ be a complex manifold with $\mbox{dim}_\C X =n$ and let $p\in\{1,\dots , n\}$. Fix a Hermitian metric $\omega$ on $X$. Assume $X$ is {\bf compact} in $(1)$ and $(2)$ below.

\vspace{1ex}

$(1)$\, A current $T\in{\cal D}^{'n-p,\,n-p}(X)$ is said to be {\bf strongly (strictly) positive} if there exists a constant $\varepsilon>0$ such that $T - \varepsilon\,\omega^{n-p}$ is strongly semi-positive as an $(n-p,\,n-p)$-current on $X$.

 We write in this case: $T>0$ (strongly) or $T - \varepsilon\,\omega^{n-p}\geq 0$ (strongly).

\vspace{1ex}

$(2)$\, A current $T\in{\cal D}^{'n-p,\,n-p}(X)$ is said to be {\bf metrically weakly (strictly) positive} if there is a constant $\varepsilon>0$ such that $T - \varepsilon\,\omega^{n-p}$ is weakly semi-positive as an $(n-p,\,n-p)$-current on $X$.

We write in this case: $T>0$ (metrically weakly) or $T - \varepsilon\,\omega^{n-p}\geq 0$ (weakly).

\vspace{1ex}

$(3)$\, A current $T\in{\cal D}^{'n-p,\,n-p}(X)$ is said to be {\bf dually weakly (strictly) positive} if \begin{eqnarray*}\int\limits_X T\wedge u >0\end{eqnarray*} for every $u\in{\cal D}_{p,\,p}(X)\setminus\{0\}$ such that $u(x)>0$ (strongly) for every $x\in\mbox{Supp}\,u$.

 We write in this case: $T>0$ (dually weakly).

\end{Def}  

Note that definitions $(1)$ and $(2)$ above can be reformulated as the following equivalences: \begin{eqnarray}\label{eqn:strong-strict-positivity_currents}T>0 \hspace{2ex} \mbox{({\bf strongly})} & \iff & \exists\,\varepsilon>0 \hspace{2ex} \mbox{such that} \hspace{2ex} \int\limits_X T\wedge u \geq \varepsilon\,\int\limits_X \omega^{n-p}\wedge u \\
  \nonumber \mbox{(on a {\bf compact} $X$)} & & \mbox{for all} \hspace{2ex} u\in{\cal D}_{p,\,p}(X) \hspace{1ex} \mbox{with} \hspace{1ex} u(x)\geq 0 \hspace{1ex} \mbox{(weakly) at every} \hspace{1ex} x\in X;\end{eqnarray}

\begin{eqnarray}\label{eqn:metrically-weak-strict-positivity_currents}T>0 \hspace{2ex} \mbox{({\bf metrically weakly})} & \iff & \exists\,\varepsilon>0 \hspace{2ex} \mbox{such that} \hspace{2ex} \int\limits_X T\wedge u \geq \varepsilon\,\int\limits_X \omega^{n-p}\wedge u \\
  \nonumber \mbox{(on a {\bf compact} $X$)} & & \mbox{for all} \hspace{2ex} u\in{\cal D}_{p,\,p}(X) \hspace{1ex} \mbox{with} \hspace{1ex} u(x)\geq 0 \hspace{1ex} \mbox{(strongly) at every} \hspace{1ex} x\in X.\end{eqnarray}

\vspace{2ex}

When $X$ is not assumed compact, the above definitions can be modified as follows: \begin{eqnarray}\label{eqn:strong-strict-positivity_currents_non-compact} \nonumber T>0 \hspace{2ex} \mbox{({\bf strongly})} & \iff & \forall\, K\subset X \hspace{1ex}\mbox{\bf compact}\hspace{2ex}  \exists\,\varepsilon_K>0 \hspace{2ex} \mbox{such that} \hspace{2ex} \int\limits_X T\wedge u \geq \varepsilon_K\,\int\limits_X \omega^{n-p}\wedge u \\
   & & \mbox{for all} \hspace{2ex} u\in{\cal D}_{p,\,p}(K) \hspace{1ex} \mbox{with} \hspace{1ex} u(x)\geq 0 \hspace{1ex} \mbox{(weakly) at every} \hspace{1ex} x\in X;\end{eqnarray}

\begin{eqnarray}\label{eqn:metrically-weak-strict-positivity_currents_non-compact}\nonumber T>0 \hspace{2ex} \mbox{({\bf metrically weakly})} & \iff & \forall\, K\subset X \hspace{1ex}\mbox{\bf compact}\hspace{2ex} \exists\,\varepsilon_K>0 \hspace{2ex} \mbox{such that}  \\
   & & \int\limits_X T\wedge u \geq \varepsilon_K\,\int\limits_X \omega^{n-p}\wedge u \\
  \nonumber & & \mbox{for all} \hspace{2ex} u\in{\cal D}_{p,\,p}(K) \hspace{1ex} \mbox{with} \hspace{1ex} u(x)\geq 0 \hspace{1ex} \mbox{(strongly) at every} \hspace{1ex} x\in X.\end{eqnarray}

Note that Observation \ref{Obs:semi-pos_same_forms-currents} implies that, whenever $T = \Omega\in C^\infty_{n-p,\,n-p}(X,\,\C)$ is a smooth form, any of the (strict) positivity properties introduced in Definition \ref{Def:strict-positivity_currents} is satisfied by $T = \Omega$ as a current if and only if it is satisfied by $T = \Omega$ at every point of $X$ as a form in the sense of $\S$\ref{subsubsection:our-def}.

\begin{Obs}\label{Obs:hierarchy_strict-pos_currents} The following implications hold for every current $T\in{\cal D}^{'n-p,\,n-p}(X)$: \begin{eqnarray*}T > 0 \hspace{2ex} \mbox{({\bf strongly})} \hspace{1ex}\implies\hspace{1ex} T > 0 \hspace{2ex} \mbox{({\bf metrically weakly})} \hspace{1ex} \implies\hspace{1ex} T > 0 \hspace{2ex} \mbox{({\bf dually weakly})}.\end{eqnarray*}

\end{Obs}  

\noindent {\it Proof.} The first implication follows at once from the definitions. To prove the second implication, fix $u\in{\cal D}_{p,\,p}(X)\setminus\{0\}$ such that $u(x)>0$ (strongly) for every $x\in\mbox{Supp}\,u$. In particular,  $u(x)\geq0$ (strongly) for every $x\in X$.

Since $K:=\mbox{Supp}\,u$ is compact, there exists $\delta_K>0$ such that $u(x) - \delta_K\,\omega^p(x)\geq 0$ (strongly) for every $x\in\mbox{Supp}\,u$. On the other hand, the metrically weak positivity hypothesis on $T$ yields an $\varepsilon_K>0$ such that the first inequality below holds: \begin{eqnarray*}\int\limits_X T\wedge u \geq \varepsilon_K\,\int\limits_X \omega^{n-p}\wedge u = \varepsilon_K\,\int\limits_K \omega^{n-p}\wedge u \geq  \varepsilon_K\delta_K\,\int\limits_X \omega^n >0.\end{eqnarray*} This completes the proof. \hfill $\Box$

\vspace{2ex}

We now obtain a duality result involving the following convex cone (denoted by ${\cal D}^{'\oplus}_{p,\,p}(X)$ in [Dem97]) that is closed in ${\cal D}^{'n-p,\,n-p}(X)$ with respect to the weak topology of currents: \begin{eqnarray}\label{eqn:strongly-semi-pos_cone_currents}\widetilde{\cal C}_X^{n-p,\,n-p}:=\bigg\{T\in{\cal D}^{'n-p,\,n-p}(X)\,\mid\, T\geq 0 \hspace{1ex} \mbox{(strongly)}\bigg\}\subset{\cal D}^{'n-p,\,n-p}(X)\end{eqnarray} in any fixed bidegree $(n-p,\,n-p)$.

\begin{The}\label{The:duality_decomp-weakly_s-semi-pos-currents} Let $X$ be a {\bf compact} complex manifold with $\mbox{dim}_\C X =n$ and let $p\in\{1,\dots , n\}$. Fix a form $\Omega\in C^\infty_{p,\,p}(X,\,\R)$ and a Hermitian metric $\omega$ on $X$. The following equivalences hold: \begin{eqnarray*}\Omega >0  \hspace{1ex} \mbox{({\bf metrically weakly}) at every point of} \hspace{2ex}  X \hspace{1ex} \iff \int\limits_X T\wedge\Omega >0 \hspace{1ex} \mbox{for every} \hspace{1ex} T\in\widetilde{\cal C}_X^{n-p,\,n-p}\setminus\{0\} \\
  \iff  \exists\,\varepsilon>0 \hspace{1ex} \mbox{such that} \hspace{1ex} \int\limits_X T\wedge\Omega\geq\varepsilon\,\int\limits_X T\wedge\omega^p \hspace{1ex} \mbox{for every} \hspace{1ex} T\in\widetilde{\cal C}_X^{n-p,\,n-p}.\end{eqnarray*}

\end{The}  

\noindent {\it Proof.} Since $X$ is compact, the first of the three conditions in the stated equivalences is equivalent to the existence of a constant $\varepsilon>0$ such that $\Omega - \varepsilon\,\omega^p \geq 0$ (weakly) at every point of $X$. By the duality (\ref{eqn:duality-currents-forms}), this is further equivalent to the third of the three conditions in the stated equivalences.

On the other hand, the third condition implies the second since $\int_X T\wedge\omega^p>0$ for every $T\in\widetilde{\cal C}_X^{n-p,\,n-p}\setminus\{0\}$. It remains to prove that the second condition implies the third.

Suppose that $\int_X T\wedge\Omega >0$ for every $T\in\widetilde{\cal C}_X^{n-p,\,n-p}\setminus\{0\}$. This amounts to the following linear continuous functional \begin{eqnarray*}l_\Omega : {\cal D}^{'n-p,\,n-p}(X)_\R\longrightarrow\R,   \hspace{5ex} l_\Omega(T) = \int\limits_X T\wedge\Omega,\end{eqnarray*} (where ${\cal D}^{'n-p,\,n-p}(X)_\R$ stands for the space of {\it real} currents of bidegree $(n-p,\,n-p)$ on $X$) being $>0$ on $\widetilde{\cal C}_X^{n-p,\,n-p}\setminus\{0\}$. Letting \begin{eqnarray*}{\cal K}_\omega:= \bigg\{T\in{\cal D}^{'n-p,\,n-p}(X)_\R\,\mid\,\int\limits_X T\wedge\omega^p = 1\bigg\},\end{eqnarray*} Proposition \ref{Prop:semi-pos-currents_b-mass_compactness} ensures that the set $\widetilde{\cal C}_X^{n-p,\,n-p}\cap{\cal K}_\omega$ is compact in ${\cal D}^{'n-p,\,n-p}(X)_\R$. We infer the existence of a constant $\varepsilon>0$ such that the restriction $(l_\Omega)_{|\widetilde{\cal C}_X^{n-p,\,n-p}\cap{\cal K}_\omega}\geq\varepsilon$. Since, for every $T\in\widetilde{\cal C}_X^{n-p,\,n-p}\setminus\{0\}$, $T/\int_XT\wedge\omega^p\in\widetilde{\cal C}_X^{n-p,\,n-p}\cap{\cal K}_\omega$, this yields: \begin{eqnarray*}\int\limits_X T\wedge\Omega\geq\varepsilon\,\int\limits_X T\wedge\omega^p \hspace{1ex} \mbox{for every} \hspace{1ex} T\in\widetilde{\cal C}_X^{n-p,\,n-p},\end{eqnarray*} as desired.  \hfill $\Box$

\vspace{2ex}

Thus, if we consider the following cone of smooth forms on $X$: \begin{eqnarray}\label{eqn:m-d_weakly--pos_cone_forms}\nonumber (m{\cal C})_X^{p,\,p} := \bigg\{\Omega\in C^\infty_{p,\,p}(X,\,\R)\,\mid\, \Omega> 0 \hspace{1ex} \mbox{(metrically weakly) at every point of}\hspace{1ex} X\bigg\}\subset  C^\infty_{p,\,p}(X,\,\R),\end{eqnarray} Theorem \ref{The:duality_decomp-weakly_s-semi-pos-currents} can be rephrased as the following 

\begin{Cor} Let $X$ be a {\bf compact} complex manifold with $\mbox{dim}_\C X =n$ and let $p\in\{1,\dots , n\}$.

\vspace{1ex}

 The cone $(m{\cal C})_X^{p,\,p}$ is the {\bf strict dual cone} of the cone $\widetilde{\cal C}_X^{n-p,\,n-p}$ under the duality (\ref{eqn:duality-currents-forms}).

\end{Cor}  

By the {\it strict dual cone} of $\widetilde{\cal C}_X^{n-p,\,n-p}$ we mean the interior (with respect to the Fr\'echet topology of $C^\infty_{p,\,p}(X,\,\R)$) of its ordinary dual cone, namely the interior of \begin{eqnarray*}\bigg(\widetilde{\cal C}_X^{n-p,\,n-p}\bigg)^\vee = \bigg\{\Omega\in C^\infty_{p,\,p}(X,\,\R)\,\mid\,\int\limits_X T\wedge\Omega\geq 0 \hspace{2ex} \forall\,T\in\widetilde{\cal C}_X^{n-p,\,n-p}\bigg\}.\end{eqnarray*}


\vspace{3ex}

\noindent {\bf References.} \\

\noindent [Buc88]\, N. P. Buchdahl --- {\it Hermitian-Einstein Connections and Stable Vector Bundles over Compact Complex Surfaces} --- Math. Ann. {\bf 280} (1988), 625-648. 

\vspace{1ex}

\noindent [Dem02]\, J.-P. Demailly --- {\it On the Frobenius Integrability of Certain Holomorphic $p$-Forms } --- In: Bauer, I., Catanese, F., Peternell, T., Kawamata, Y., Siu, YT. (eds) Complex Geometry. Springer, Berlin, Heidelberg. \url{https://doi.org/10.1007/978-3-642-56202-0_6}

\vspace{1ex}

\noindent [Dem97]\, J.-P. Demailly --- {\it Complex Analytic and Algebraic Geometry} --- \url{https://www-fourier.univ-grenoble-alpes.fr/~demailly/manuscripts/agbook.pdf}

\vspace{1ex}

\noindent [DP25a]\, S. Dinew, D. Popovici --- {\it $m$-Positivity and Regularisation} --- arXiv:2510.25639v1 [math.DG].

\vspace{1ex}

\noindent [DP25b]\, S. Dinew, D. Popovici --- {\it $m$-Pseudo-effectivity and a Monge-Amp\`ere-Type Equation for Forms of Positive Degree} --- arXiv:2510.27362v1 [math.DG].

\vspace{1ex}

\noindent [Kob82]\, S. Kobayashi --- {\it Curvature and Stability of Vector Bundles} --- Proc. Jap. Acad. {\bf 58} (1982), 158-162. 

\vspace{1ex}

\noindent [Kob87]\, S. Kobayashi --- {\it Differential Geometry of Complex Vector Bundles} --- Princeton University Press, 1987.

\vspace{1ex}

\noindent [LY86]\, J. Li, S-T Yau --- {\it Hermitian-Yang-Mills Connection on Non-K\"ahler Manifolds} --- Mathematical aspects of string theory (San Diego, Calif., 1986), 560–573, Adv. Ser. Math. Phys., {\bf 1}, World Sci. Publishing, Singapore, 1987.

\vspace{1ex}

\noindent [Lub83]\, M. L\"ubke --- {\it Stability of Einstein-Hermitian Vector Bundles} --- Manuscripta Math. {\bf 42} (1983), 245-257.

\vspace{1ex}

\noindent [LT95]\, M. L\"ubke, A. Teleman --- {\it The Kobayashi-Hitchin Correspondence} --- World Scientific, 1995.

\vspace{1ex}

\noindent [Mic83]\, M. L. Michelsohn --- {\it On the Existence of Special Metrics in Complex Geometry} --- Acta Math. {\bf 143} (1983) 261-295.

\vspace{1ex}

\noindent [Xia26]\, M. Xia --- {\it Transcendental $b$-Divisors II --- Monotonicity Theorem} --- arXiv:2603.14362v1

\vspace{6ex}

\noindent Institut de Math\'ematiques de Toulouse,

\noindent Universit\'e de Toulouse,

\noindent 118 route de Narbonne, 31062 Toulouse, France

\noindent Email: popovici@math.univ-toulouse.fr

\end{document}